\documentclass[12pt,a4paper]{article}

\usepackage{pdflscape}
\usepackage{amsmath}
\usepackage{amssymb}
\usepackage{latexsym}
\usepackage{srcltx}
\usepackage{graphics}
\usepackage{color}
\usepackage{epsfig}
\usepackage{color}
\usepackage{subcaption}
\usepackage[export]{adjustbox}
\usepackage{multirow}

\newcommand{\R}{{\mathbb R}}
\newcommand{\re}{{\mathbb R}}
\newcommand{\n}{{\mathbb N}}

\newcommand{\cA}{{\mathcal{A}}}

\newcommand{\cV}{{\mathcal{V}}}

\newcommand{\cT}{{\mathcal{T}}}

\newcommand{\cP}{{\mathcal{P}}}

\newcommand{\bx}{{\boldsymbol{x}}}
\newcommand{\by}{{\boldsymbol{y}}}

\newcommand{\bp}{{\boldsymbol{p}}}
\newcommand{\bq}{{\boldsymbol{q}}}
\newcommand{\ba}{{\boldsymbol{a}}}
\newcommand{\bb}{{\boldsymbol{b}}}

\newcommand{\bu}{{\boldsymbol{u}}}

\newcommand{\bh}{{\boldsymbol{h}}}

\newcommand{\nill}{{\boldsymbol{O}}}
\newcommand{\bell}{{\boldsymbol{\ell}}}

\newcommand{\vardot}{\mathord{\,\cdot\,}}

\newtheorem{theorem}{Theorem}
\newtheorem{prop}{Proposition}
\newtheorem{lemma}{Lemma}

\newtheorem{remark}{Remark}
\newtheorem{ex}{Example}
\newtheorem{defi}{Definition}

\date{}

\date{}

\author{Vladimir Protasov 
	\thanks{
			Moscow  State University, {e-mail: \tt\small
			v-ptotassov@yandex.ru}} , 
	Rinat Kamalov
	\thanks{Moscow Institute of Physics and 
		Technology, Russia; Gran Sasso Science Institute, Italy,  {e-mail: \tt\small
			rinat020398god@yandex.ru}} 
}

\title{Chebyshev approximation by non-Chebyshev systems}

\begin{document}

\maketitle

\begin{abstract}
We address the problem  of the best  uniform approximation 
by linear combinations  of a finite system of functions. 
If the system is  Chebyshev and the problem is unconstrained, then 
the classical Remez algorithm provides a fast and precise solution. 
For non-Chebyshev systems, this problem may offer a great resistance.  The same 
happens to approximations under linear constraints. 
We propose a solution 
by modifying the concept of  alternance  and of the Remez iterative procedure. 
A criterion of the best approximation is proved and  
the full set of polynomials of best approximation 
(which may not be unique in the non-Chebyshev case) 
is characterized. The method of finding the best polynomial 
is applicable for arbitrary functional systems under arbitrary linear constraints. 
The efficiency is demonstrated in examples with 
 systems of complex exponents, Gaussian functions, and  lacunar  polynomials. 
 As an application, the Markov-Bernstein type inequalities are obtained 
for those systems.
Applications to signal processing, linear ODEs, switching dynamical systems are considered. 

\bigskip


\smallskip

\smallskip

\begin{flushright}
\end{flushright}

\end{abstract}

\section{Introduction}
The uniform finite-dimensional approximation of functions on compact domains is well-understood  for approximations by Chebyshev systems on a segment. 
In this case, for every~$f\in C[a,b]$, the closest to~$f$
polynomial is always unique and is characterized by the alternance criterion (the Chebyshev - Vall\`ee-Poussin  theorem). 
Finding  the polynomial of best approximation is realized by the Remez algorithm~\cite{R}. 
A system of~$n$ functions on a set~$K$ is {\em Chebyshev} (also called {\em Haar}) if  
every nontrivial polynomial by this systems, i.e., a linear combination of its functions, 
has at most~$n-1$ zeros on~$K$. 
They are, for example,  
algebraic polynomials on a segment and trigonometric polynomials
on the period. See~\cite{DS, KS, KN} for the general theory of Chebyshev approximations. 

For non-Chebyshev systems, the situation is 
different. Neither the alternance criterion nor the Remez algorithm works. 
Very little is known on the  uniform polynomial approximation in this case, especially from a numerical aspect. 
Many important problems in signal processing, numerical ODEs, etc. 
deal with non-Chebyshev systems such as quasipolynomials, band-limited functions, 
lacunary algebraic  polynomials (polynomials spanned by 
some integer powers~$t^{m_k}, \, k=1, \ldots , n$), M\"untz polynomials
(spanned by 
arbitrary nonnegative  powers), etc. All multivariate systems  are, as a rule, non-Chebyshev. 
Another class of problems which is not covered by the Chebyshev analysis 
is the constrained uniform approximation. For example, an algebraic  polynomial of best approximation with a prescribed value of the derivative at some point may not possess an alternance 
and may not be unique. Actually, imposing a linear constraint can destroy all properties of 
the Chebyshev approximation. 

In this paper we begin with a criterion for the best approximation by polynomials 
spanned by an arbitrary system of continuous functions on a compact domain. 
It is formulated with the concept of   ``generalized alternance''. 
This  result is actually not new and appeared in the literature in different forms  (we give references in the next section). 
Then we prove the classification theorem that characterize the whole set 
of polynomials  of best approximations in the case of non-uniqueness. 
Then we address a constrained approximation problem: formulate and prove  the criterion and the classification for uniform approximation under  linear constraints. 
Using these results we derive a ``Remez-like'' method for numeical computation of the 
best approximating polynomials. The main difficulty is that on some steps the algorithm may 
arrive at  degenerate configurations, in which case it slows down or even gets stuck.
Besides, 
the alternance can be degenerate, which leads to non-uniqueness  of the 
solution and to difficulties in the numerical implementation  (linear systems solved in each iteration become ill-conditioned, etc.) We present a method
of regularization that overcomes both those problems. The performance of 
the new approach is shown through numerical experiments. In the vast majority of  examples, including degenerate ones,  the polynomial of best approximation is found within a few iterations. This is also demonstrated in applications from signal processing and from the stability of dynamical systems. 


\bigskip

\section{The roadmap of the main results} We consider an arbitrary metric compact set~$K$ and a finite system of continuous linearly independent  real-valued functions~$\Phi = \{\varphi_1, \ldots , \varphi_n\}$ on it. 
A linear combination~$p = \sum_{k=1}^n p^{\,k} \, \varphi_k \, , \ p^{\,k}\in \re$, will be referred to as 
a {\em polynomial by the system~$\Phi$}. All polynomials form a linear space~$\cP$, which is 
naturally identified with~$\re^n$. The vector of coefficients of the polynomial~$p$
will be denoted by the bold letter:~$\bp = (p^1, \ldots , p^n) \in \re^n$. 
All the  results 
are also valid for unbounded domains provided that all~$\varphi_k(t)$ tend to zero as~$t\to \infty$.

By~$C(K)$ we denote the space of continuous functions on a metric compact set~$K$ with a standard uniform norm. For a function~$f\in C(K)$, 
the problem of the {\em best uniform polynomial approximation} is
\begin{equation}\label{eq.mail}
	\left\{
	\begin{array}{l}
		\|p -f \|_{C(K)} \ \to \ \min\\
		p\ \in \ \cP. 
	\end{array}
	\right.
\end{equation}
Note that the objective function is continuous and coercive on~$\cP$, hence, the optimal 
polynomial~$\hat p$ always exists, although may not be unique.  
The problem~(\ref{eq.mail}) is convex,  however, the standard tools of convex optimization
are usually less efficient  than special approximating algorithms. 
We will discuss this issue in Section~10. 
If~$\Phi$ is a Chebyshev system on a segment~$[a,b]$, then $\hat p$ is unique and 
can be found  by  the Remez algorithm, which has a linear convergence and normally computes the 
solution fast and with a high accuracy. It is based on the 
classical Chebyshev criterion: $\hat p$ is a polynomial 
of the best approximation precisely when there exists a sequence of~$n+1$ points 
$\tau_1 < \ldots < \tau_{n+1}$ called {\em alternance} such that the difference~$\hat p(t) - f(t)$ reaches the maximal absolute value at the points~$t=\tau_j$ and its sign alternates 
in~$j$, i.e., $\hat p(\tau_j) - f(\tau_j) \, = \, (-1)^{j+\delta}\|\hat p - f\|$, 
where~$\delta$ is either zero or one. Let us briefly describe the Remez algorithm. 
By the  $k$th iteration it produces points~$\bigl\{t^{(k)}_i\bigl\}_{i=1}^{n+1}$
and a polynomial~$p_k$ such that the values~$p_k\bigl(t_i^{(k)}\bigr) - f\bigl(t_i^{(k)}\bigr)$
are equal in modulus (not necessary equal to~$\|p_k-f\|$) 
and their signs alternate in~$i$. We take  the point of maximum of~$|p_k-f|$, 
denote it by~$t_0^{(k)}$, and   replace some point~$t^{(k)}_s$ by~$t_0^{(k)}$. 
Then we construct the next polynomial~$p_{k+1}$ which possesses the alternating property on the new system of points. The point~$t^{(k)}_s$ is chosen to be the neighboring node to~$t_0^{(k)}$ with the same sign of the difference~$p-f$. 

The Chebyshev property holds for several important functional systems such as 
algebraic and trigonometric polynomials, real exponents, etc. However, in general, Chebyshev systems are rather  exceptional. For arbitrary systems, neither of the aforementioned facts take place.
The optimal polynomial~$\hat p$ is not necessarily  unique and may  not possess an alternance. 
The ``Refinement theorem'' (Theorem~B in Section~3) implies the  existence of  
a set of points~$\{\tau_i\}_{i=1}^N, \, N\le n+1$, such that 
no polynomial from~$\cP$ has all values~$|p(\tau_i)-f(\tau_i)|$ smaller than~$\|\hat p - f\|$. 
This set can be considered as  a weakened version of alternance,   however, it may be degenerate
(when~$N < n+1$) and the signs of~$p(\tau_i)-f(\tau_i)\, , \, i=1, \ldots , n+1$, 
may not alternaite. Therefore, this ``weakened alternance'' is not sufficient 
to characterize the best approximation. Moreover, the routine of Remez' algorithm becomes inapplicable, because the choice of the 
point~$t_s^{(k)}$ to be replaced by~$t_0^{(k)}$ is not valid without the sign alternating.   

Theorem~\ref{th.10} proved in Section~3 formulates the best
approximation criterion in the case of a non-Chebyshev system.   The signs alternating  property is replaced by the following  geometric condition: 
the simplex formed by  the  oriented moment vectors in~$\re^n$ contains the origin. 
Here we need some further notation.   
For an arbitrary point~$t\in K$, we consider the {\em moment vector}  $\bu(t)\, = \, \bigl(\varphi_1(t), \ldots ,  \varphi_n(t)\bigr)^T$. If, in addition, we are given a 
function~$f\in C(K)$ and a polynomial~$p\in \cP$, then 
we define the {\em oriented moment vector}
$\ba(t) \, = \, \sigma (t)\, \bu(t)$, where~$\sigma(t) \, = \ {\rm sign}\, \bigl(p(t) - f(t)\bigr)$. For a given polynomial~$\hat p \in \cP$, the {\em generalized alternance} is a set~$\{\tau_i\}_{i=1}^m \subset K, \, 
m\le n+1$, 
such that~$\hat p(\tau_i) - f(\tau_i)\, = \, \sigma_i\, \|\hat p-f\|\, , \
\sigma_i \in \{-1, 1\},$ for all~$i$, and the $(m-1)$-dimensional  simplex~$\Delta$ with vertices~$\{\sigma_i\bu (\tau_i)\}_{i=1}^m$ contains the origin in~$\re^n$. In other words, some convex combination of the vectors 
$\ba(\tau_i), \, i=1, \ldots , m,$ is equal to zero.  

Theorem~\ref{th.10} asserts that 
$\hat p$ is a best approximation polynomial if and only if it possesses 
the generalized alternance. The results similar to Theorem~\ref{th.10} 
appeared in the literature in various forms, see~\cite[p.12]{CL}, \cite[theorem~5.2]{CK},~\cite[theorem 4.2]{SS}, \cite{SU}.  
Its proof follows in a rather 
straightforward manner by applying the Refinement theorem and the Farkas lemma (see Section~1).
Another possible proof  involves subdifferentials (Remark~\ref{r.20}). 

If the system is non-Chebyshev, then the best approximation may not be unique
as the following simple example demonstrates: 
\begin{ex}\label{ex.10}
	{\em Let~$K = [-1,1], \, \Phi = \{t, t^2, \ldots , t^n\}$. 
		For the function~$f=1$, each of the functions~$t^k$ is a best approximation polynomial. 
		Indeed, $p(0) = 0$ for all~$p\in \cP $, hence, the distance from~$f$ to $\cP$ is at least one. 
		On the other hand, it is equal to one for every~$p_k(t) = t^k$. It is interesting that for each~$p_k$, the generalized alternance is a singleton~$\{0\}$. On the other hand, if~$n$ is even 
		and~$T_n = \cos \, (n \, {\rm acos}\, t)$ is the Chebyshev polynomial, then 
		$p(t) = 1 - (-1)^{n/2}T_n$ is also a best approximation polynomial, 
		and its alternance~$\tau_k = \cos \frac{\pi k}{n}, \, k=0, \ldots , n$, is the same as for the Chebyshev polynomial~(Example~\ref{ex.40} in Section~3).  
		
	}
\end{ex}

If, for a given function~$f$,  the best approximation polynomial is not unique, 
how can we find them all? A complete classification is given in Theorem~\ref{th.15}. 
The classification  is surprisingly simple: all polynomials of the best approximation must have the same generalized alternance with the same values on it. Note that each polynomial 
may also have other alternances. 

Theorems~\ref{th.10} and~\ref{th.15} make it possible to generalize Remez' algorithm 
for an arbitrary system of functions~$\Phi$ on a  compact set. In each iteration we 
take the point~$t_0$ of maximum modulus~$|p_k(t) - f(t)|$ and replace
one of the vertices of the simplex~$\Delta_k$ with $\sigma_0\bu(t_0), \, 
\sigma_0 = {\rm sign}\, \bigl(p_k(t_0) - f(t_0)\bigl),$ so that the new simplex 
$\Delta_{k+1}$ still contains the origin. The algorithm halts when 
we localize the distance from~$f$ to~$\cP$ to 
a segment of a given length~$\varepsilon$.  This algorithm has a 
linear rate of convergence in the nondegenerate case. The latter condition is  significant: 
in many situations 
the simplex~$\Delta_k$ that appears at the~$k$th iteration may be close to degenerate, 
after which the algorithm slows down or stalls. In this case algorithm also 
faces technical problems because 
the linear systems solved in each iteration become ill-conditioned. Moreover, if~$m < n+1$, then  this simplex always degenerates for large~$k$
since the limit simplex is~$(m-1)$-dimensional. To resolve these issues 
we present 
the method of regularization  and modify the algorithm.
The key novelty is the choice 
of the new point~$t_0^{(k)}$, see~Section~5. The efficiency of the new method 
is confirmed numerically (Section~9).    
\smallskip 

In Section~6  we address the constrained  approximation. 
For   given linear functionals~$\ell_j: \cP \to \re$ and numbers~$b_j\, , \ j=1, \ldots , r$, we 
consider the problem 
\begin{equation}\label{eq.constr}
	\left\{
	\begin{array}{l}
		\|p -f \|_{C(K)} \ \to \ \min\\
		\ell_j(p) \ = \ b_j, \ j = 1, \ldots , r,\\ 
		p\, \in \, \cP. 
	\end{array}
	\right.
\end{equation}
Even if the system~$\{\varphi_i\}_{i=1}^n$ is Chebyshev, problem~(\ref{eq.constr}) 
cannot be solved by the 
standard approach: the optimal polynomial may not have an alternance (even generalized one)
and the Remez algorithm may not be applicable.  
\begin{ex}\label{ex.15}
	{\em Let~$\cP$ be the space of algebraic polynomials of degree~$n$ on the segment~$[-1,1]$. 
		The problem of approximating the identical zero~$f\equiv 0$ 
		under the constraint $\ell(p) = p_0 = -1$ (the constant coefficient of~$p$ is equal to~$-1$) 
		is equivalent to the problem from Example~\ref{ex.10}.  Hence,  
		each of the functions~$t^k, \, k\ge 1$,  is a best approximation polynomial with the one-point 
		alternance~$\{0\}$. We see that the Chebyshev system~$\{t^k\}_{k=0}^n$ loses its properties after imposing a linear constraint.  
		
	}
\end{ex}
Nevertheless, it turns out that the generalized alternance approach and the 
regularization technique can be adopted to constrained problems and can efficiently solve them. 
Theorems~\ref{th.10c} and \ref{th.15c} characterize polynomial of the 
best approximation. Then we present the algorithm that can be used to compute them.   

Section~8 deals with several applications. Particular attention is paid to the signal processing, where signals often have the form of 
quasipolynomials 
(real parts of complex exponentials). If the signal is given by observation, with possible errors and noise, 
the problem is to recover it in an optimal (in some sense) way. 
Our method restores it with the minimal uniform distance from the observed function. 
If, in addition, the signal is given by means of a linear sampling, 
then the optimal recovery  is done by solving the constrained problem~(\ref{eq.constr}). 
The same technique works for other types of signals, for instance, 
the ones spanned  by shifts and contractions of the Gaussian function or of the Cauchy 
function. 

Other applications are on the stability of linear dynamical systems and analyzing the growth of their trajectories. 
\bigskip

\section{The generalized alternance} We begin with a criterion of the best uniform approximation. Let~$K$ be a metric compact space, $f\in C(K)$, $\Phi = \{\varphi_i\}_{i=1}^n$ be a  system of linearly independent functions from~$C(K)$
and $\cP$ be the space of polynomials by the system~$\Phi$. We denote with~$\, {\rm co}\, (X)$
we denote the convex hull of a set~$X$. 
\begin{defi}\label{d.10}
	For given~$p\in \cP$,  a set of points~$\{\tau_1, \ldots , \tau_{m}\}, \, m\le n+1$, 
	is called {\em generalized alternance} if 
	\smallskip 
	
	1) $\quad p(\tau_i) - f(\tau_i)\ = \ \sigma_i \|p-f\|,  \  
	\sigma_i \in \{-1, 1\}, \quad $ for all~$i =1, \ldots , m$,
	;  
	\smallskip 
	
	2) $\ $ the convex hull of the vectors~$\{\sigma_i\bu(\tau_i)\}_{i=1}^m$ in~$\re^n$, where \\ $ \bu = \left(\varphi_1(t), \ldots, \varphi_n(t)\right)^{T} $,
	contains the origin. 
\end{defi}
\medskip 

Using the notation~$\ba_i = \sigma_i \bu(\tau_i)$, we have  
$\nill \, \in \, \Delta  \, = \, {\rm co}\, \{\ba_1, \ldots , \ba_m\}$. 
In what follows we deal with generalized alternance only and  drop the word ``generalized''. 

\begin{theorem}\label{th.10}
	A polynomial~$p \in \cP$ is of the best approximation for a function~$f \in C(K)$ if and only if 
	it has an alternance. 
\end{theorem}
As we know, for a Chebyshev system on a segment the alternance possesses the sign alternating property. 
For general systems this property is replaced by that~$\,\nill \in \Delta$. In fact, already for 
two-element systems~$\Phi$ on~$[-1,1]$, there are examples of polynomials of the best approximation
with non-alternating signs of~$p(\tau_i) - f(\tau_i)$.
Example~\ref{ex.50} at the end of this section demonstrates this situation 
for the system~$\Phi = \{t, t^2\}$. 

The proof of Theorem~\ref{th.10} is based on two 
results of convex analysis formulated below. The first one is the well-known  Farkas lemma, see, for instance~\cite{Roc}. 
\smallskip 

\noindent \textbf{Lemma~A} {\em (Farkas' lemma)}
{\em For an arbitrary set of vectors~$\bb_1, \ldots , \bb_N$ in~$\re^N$
	the following are equivalent: 
	\smallskip 
	
	1) their convex hull contains the origin; 
	\smallskip 
	
	2) the system of linear inequalities~$(\bb_i, \bx) < 0, \ i=1, \ldots , N$, does not have a solution
}
\smallskip 

\noindent The second auxiliary result is the 
Refinement theorem; its proof can be found in~\cite{MIT}.   
\smallskip 

\noindent \textbf{Theorem~B} {\em (Refinement theorem)}
{\em Let  $K$ be a compact subset of a metric space, 
	a function $\, F: K \times \re^n\, \, \to \, \re$ be such that  
	for each~$t\in K$, the function $F(t, \vardot )$ is convex and 
	for each $\bp \in \re^n$, the function  
	$F(\vardot , \bp)$ is continuous. Then there exist~$N \le n+1$ and 
	points $\tau_1, \ldots , \tau_N \in K$ such that 
	$$
	\inf_{\bp\in \re^n} \, \left(\, \max_{t \in K} \, F(t,\bp)\, \right)\ = \
	\inf_{\bp\in \re^n} \, \left(\, \max_{i=1, \ldots , N} \, F(\tau_i,\bp)\, \right)\, .
	$$
}
\smallskip

Applying Theorem~B to the function~$F(t, \bp) \, = \, |p(t) - f(t)|$, 
where~$p$ is the polynomial from~$\cP$ with the vector of coefficients~$\bp$, we see that 
finding the best approximation polynomial on~$K$ can always be reduced to that on some finite subset $\{\tau_i\}_{i=1}^m \subset K$. Therefore, some kind of ``weakened alternance'' always exists. 
The problem is that Theorem~B asserts only the existence and does not give any recipe how to find such points. Nevertheless, it gives a bridge to Theorem~\ref{th.10} 
for charactering polynomials of the best approximation. 

\medskip

{\tt Proof of Theorem~\ref{th.10}}. {\em Sufficiency}. If a polynomial~$p$ possesses an   alternance $\{\tau_j\}_{j=1}^m$ but does not provide the best approximation, then 
there exists~$h\in \cP$ such that~$\|p+h - f\| \, < \, \|p - f\|$. 
Therefore, at each point~$\tau_i$, we 
have
$$
\bigl|\, p(\tau_i)- f(\tau_i) +h(\tau_i) \bigr| \ < \ \bigl\|p - f\bigr\| \ = \ \bigl|p(\tau_i) - f(\tau_i)\bigr| \, .
$$ 
This implies that the value~$h(\tau_i)$ has an opposite sign to~$(p - f)(\tau_i) \, = \, 
\sigma_i \, \|p-f\|$. Hence, ${\rm sign}\, h(\tau_i) \, = \, - \, \sigma_i$ and 
consequently,~$\, \sigma_ih(\tau_i) \, < \, 0$. If $\bh$ denotes the vector of coefficients
of the polynomial~$h$, then~$h(\tau_i) \, = \, \bigl(\bh, \bu(\tau_i) \bigr)$ and 
$\, \sigma_i\, h(\tau_i) \, = \, \bigl(\sigma_i \bh,  \bu(\tau_i) \bigr)\, = \, 
\bigl(\bh, \sigma_i \bu(\tau_i) \bigr) \, = \, \bigl(\bh, \ba_i \bigr)$. We see that~$\bigl(\bh, \ba_i \bigr) < 0$ for all~$i = 1, \ldots , m$, which 
contradicts to Farkas' lemma. 
\smallskip 

{\em Necessity}. Let $\hat p$ be the best approximation polynomial. Applying 
the Refinement theorem to the function~$F(t, \bp) \, = \, |p(t) - f(t)|$, we obtain a set of 
points~$\tau_i \in K, \, i = 1, \ldots , N$, where $N\le n+1$, 
such that~$\hat p$ gives the best approximation of the function~$f$ 
on this set. After renumbering it can be assumed that
$|\hat p(\tau_i) -f(\tau_i)| = \|\hat p - f\|$ for~$i \le m$ and
$|\hat p(\tau_i) -f(\tau_i)| < \|\hat p - f\|$ for~$i > m$, where 
$0 \le m\le N$. If $m=0$, then~$\max_{i=1}^N
|\hat p(\tau_i) -f(\tau_i)| < \|\hat p - f\|$, which is impossible, 
hence~$m\ge 1$.    
If there exists a polynomial~$h \in \cP$
such that for each~$i=1, \ldots , m$, the value~$h(\tau_i)$ has the sign~$ \, -  \, \sigma_i$,
then for sufficiently small~$\lambda > 0$, we have
$|\hat p(\tau_i)  - f(\tau_i) + \lambda h(\tau_i)| \, <  \, |\hat p(\tau_i) - f(\tau_i)|$
for all~$i=1, \ldots , N$, which 
contradicts to the best approximation property of~$\hat p$ on the set~$\{\tau_i\}_{i=1}^N$.  
Thus, there is no polynomial~$h\in \cP$ for which~$\sigma_i h(\tau_i) < 0, \ 
i=1, \ldots , m$. In terms of the corresponding vectors of coefficients~$\bh, \bp$, we have 
$\, \bigl(\bh , \ba_i \bigl)\, =  \, 
\bigl(\bh , \sigma_i\, \bu(\tau_i) \bigl) \, = \, \sigma_i \, h(\tau_i)\, < \, 0$. 
Hence, the system of linear 
inequalities~$\, \bigl(\bh , \ba_i \bigl)\,  < \, 0, \, i = 1, \ldots , m$, 
has no solution. 
By the Farkas lemma, $\nill \in {\rm co}\, \{\ba_1, \ldots , \ba_n\}$.

{\hfill $\Box$}
\medskip 

\begin{remark}\label{r.20}
	{\em Another proof of Theorem~\ref{th.10} can be done by involving  subdifferentials. 
		Let us describe the idea. The point~$\hat p \in \cP$ gives an absolute minimum of the 
		function $F(p)\, = \, \|p-f\|$ if and only if~$\, \nill  \, \in \, \partial F (\hat p)$, where
		$\partial F$ denotes the subdifferential.  
		Since~$F(\hat p) = \max_{t\in K}\, |\hat p(t) - f(t)|$, it follows that 
		the subdifferential~$\partial F(p)$ can be computed by the Dobovitsky-Milyutin theorem 
		on the suddifferential of maximum. Denote
		$$
		\Gamma \ = \ \Bigl\{ \, 
		\sigma(t)\bu(t) \, \in \re^n\ : \ t\in K, \ |\hat p(t) - f(t)| \, = \, \|\hat p -f\|\, \Bigr\}.
		$$ 
		We obtain~$\partial F (\hat p)\, = \, {\rm co}\, \Gamma$. 
		Thus, $\nill \in {\rm co}\, \Gamma$. Now we invoke  the Caratheodory theorem 
		and conclude that there exist 
		$m\le n+1$  points from~$\Gamma$, whose convex hull contains~$\nill$.  
	}
\end{remark}

Clearly, for every~$f\in C(K)$, the set of polynomials of the best approximation is a convex compact subset of~$\cP$. Due to Theorem~\ref{th.10}, each element~$p$ of this subset 
possesses an alternance, which is not necessarily unique and depends on~$p$. 
It turns out, however,  that all those alternances have a nonempty intersection, which is a common alternance for all polynomials of the best approximation.  
\begin{theorem}\label{th.15}
	Let~$\cP$ be the space of polynomials by a system~$\Phi$. 
	For every~$f\in C(K)$, all polynomials of the best approximation have the same alternance with the same values on it.  
\end{theorem}
This means that for every function~$f$, there exists a set 
of points~$\{\tau_i\}_{i=1}^m$ such that all polynomials  
of the best approximation have this set as an alternance and 
for every~$i$, the value of all those polynomials at the point~$\tau_i$
is the same. 
\smallskip 

{\tt Proof}. Let~$\, \{\, p_{v} \in \cP \, : \ v\in \cV\, \}$ be the set of all 
polynomials of the best approximation,~$\cV$ is some index set. Without loss of generality it can be assumed that~$\|p_v - f\| = 1,\ v\in \cV$. 
Every point of absolute maximum of~$|p_v - f|$ will be called~{\em extreme point}
of~$p_v$.
Thus,  at every extreme point~$\tau$, we have $\, 
p_v(\tau) - f(\tau) \, = \, \sigma_v(\tau) \in \{-1, 1\}$. 
Denote by~$\cT_{v}$ the set of pairs
~$\bigl\{ \, (\tau, \sigma_v(\tau))\, : \ \tau \, \mbox{is an extreme point of}\, p_v\bigr\}$.  
By Theorem~\ref{th.10}, for every~$v \in \cV$, the set 
of extreme points of~$p_v$ contains an alternance of~$p_v$, therefore, 
the set~$C_v = {\rm co}\, \{\ba(\tau) : \, \tau \, \mbox{is an extreme point of}\, p_v\}$ 
contains the origin. Since all polynomials of the best approximation 
form a convex compact subset in~$\cP$, for every~$v_1, \ldots , v_N \in \cV$, the polynomial~$\frac{1}{N}\, \sum_{j=1}^N 
p_{v_j}$ belongs to this set. Denote it by~$p_a, \, a\in \cV$
and use the short notation~$v_j = j$. 
The convexity implies that for every extreme point~$\tau$ of~$p_a$, 
all the values 
$p_j(\tau) - f(\tau), \, j=1, \ldots , j$ are equal~(either to~$1$ or to~$-1$), 
hence, $\tau$ is an extreme point for all~$p_j$ with the same sign~$\sigma_j(\tau) = \sigma_a(\tau)$. Therefore,  
$\cT_v  \in \bigcap_{j=1}^N\, \cT_{j}$. We see that  the family of compact 
sets~$M = \{\cT_{v}\, : \, v\in \cV\}$ possesses the following property: 
the intersection of an arbitrary finite subfamily of~$M$ contains the set from~$M$. 
Then  all~$\cT_v \in M$  have a nonempty intersection, which  is 
also a set~from~$M$. Let it be~$\cT_{b}, \, b\in \cV$. 
The set~$C_v$ contains the origin, therefore, by Caratheodory's theorem~\cite{MIT}, 
there are extreme points~$\tau_1, \ldots , \tau_m$ of~$p_b\, , \ m\ne n+1$
such that~$\nill \in {\rm co}\, \{\, \sigma(\tau_i)\bu(\tau_i)\, \}_{i=1}^m$.  Since they are extreme points 
for every~$p_v$ with the same signs~$\sigma_v(\tau_i), \, v \in \cV$, 
we see that they form an alternance for~$p_v$.

{\hfill $\Box$}
\medskip 

\begin{ex}\label{ex.40}(A one-point alternance). 
	{\em In Example~\ref{ex.10}, all polynomials~$p(t) = \sum_{k=1}^n p^kt^k$ of the best approximation for~$f\equiv 1$ on the segment~$[-1,1]$ have the same alternance~$\{0\}$. 
		On the other hand, some of them have other alternances. 
		For instance, the polynomial $p(t) = 2t^2$ also gives the best approximation and 
		it has  an alternance~$\{-1, 0, 1\}$; the polynomial   
		$p(t) = 1 - (-1)^{n/2}T_n$, where~$n$ is an even number and~$T_n$
		is the Chebyshev polynomial, have the alternance~$\cos \, \frac{\pi k}{n}, \, k = 
		0, \ldots , n$ (the same as~$T_n$ has). 
	}
\end{ex}
\begin{ex}\label{ex.50}(An  alternance with non-alternating signs). 
	{\em  Consider the space of quadratic polynomials without the constant term:~$\Phi \, = \, \{\varphi_1, \varphi_2\}\, = \, \{t^2, t\}$
		and approximate the function~$f(t) = t^4 + t^3 - \frac14$ on the segment~$[-1, 1]$. 
		Let us show that~$p(t) = \frac34 t^2 + \frac12 t$ is a unique polynomial 
		of the best approximation. We have 
		$$
		p(t) \, - \, f(t)\ = \ \frac12 \ - \ \Bigl(\, t+1\, \Bigr)^2\, \Bigl(\, t-\frac12\, \Bigr)^2\, . 
		$$
		The modulus of this function reaches its maximum on~$[-1,1]$ at three points: $\tau_1 = -1, \tau_2 = \frac12$  (at those points~$\, p  -  f = \frac12$), and $\tau_3 = 1$, where~$p-f = -\frac12$, see Fig.~\ref{fig-nosign}.

		\begin{figure}[h]
			\caption{}
			\centering
			
			\begin{subfigure}{0.4\textwidth}
				\centering
				\includegraphics[width=\linewidth]{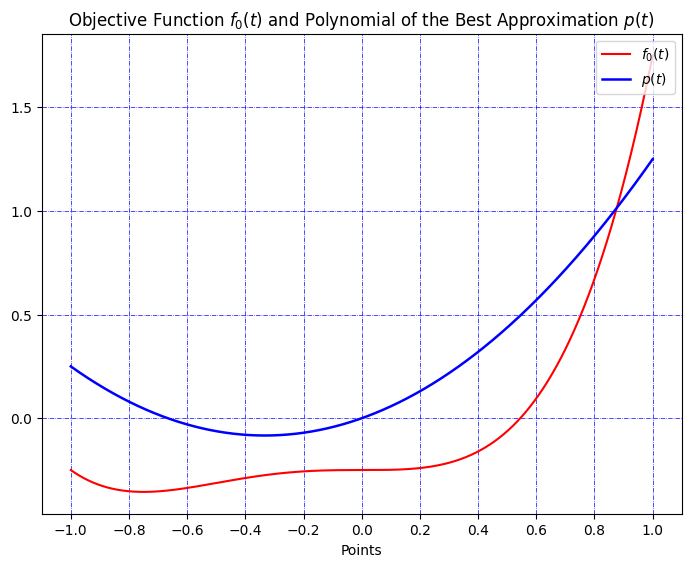}
				\caption{$ f(t) $ and $ p(t) $}
			\end{subfigure}\hfill
			\begin{subfigure}{0.4\textwidth}
				\centering
				\includegraphics[width=\linewidth]{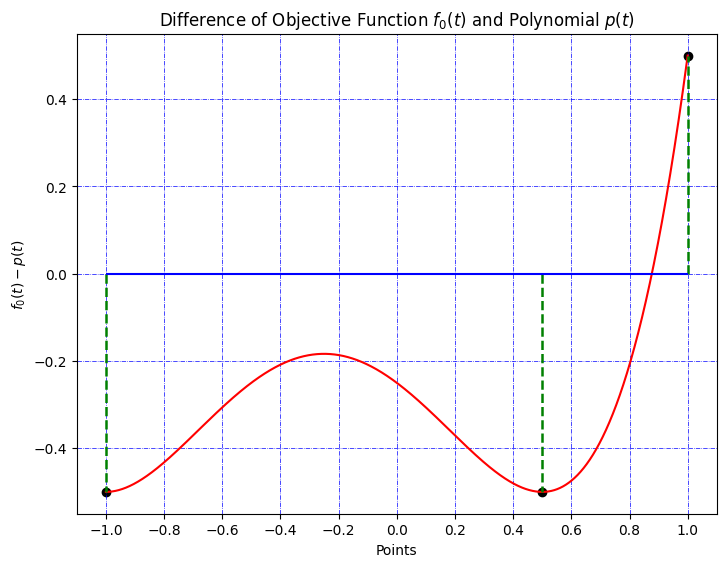}
				\caption{$ f(t) - p(t) $}
			\end{subfigure}\hfill
			
			\label{fig-nosign}
		\end{figure}
		
		This is an alternance because~$\ba_1 = \bigl(-1,1\bigr)^T; 
		\ba_2 = \bigl(\frac12 ,\frac14\bigr)^T, \, \ba_3 = - (1,1)^T$ and 
		$\nill \, \in {\rm co} \, \{\ba_1, \ba_2 , \ba_3\}$. Thus,~$p$ is a polynomial 
		of the best approximation. It is unique by Theorem~\ref{th.15}.  Indeed, $p$ has a unique alternance, since~$\tau_i$ are the only points of 
		maximum of the function~$|p-f|$ and there is no subset of those points for which 
		the convex hull if the corresponding vectors~$\ba_i$ contains the origin. Consequently, another best approximation polynomial, if it exists, has the same alternance~ and the same values on it, which is impossible since
		there is only one quadratic function passing through three fixed points.  
		Note that the difference~$f-p$ has the same sign at~$\tau_1$ and~$\tau_2$.  
		
	}
\end{ex}
%
\bigskip

\section{Algorithm for the best approximation. The regular case}

We derive Algorithm~1 that follows the Remez scheme, but with another 
rule of choosing the next point of alternance at the~$k$th  iteration. 
The sign alternating property is replaced by the requirement that~$\nill \in \Delta_k$, 
where  
the simplex $\Delta_k$ has vertices at the oriented moment vectors~$\ba_i \, = \, 
\sigma(\tau_i)\bu(\tau_i)$. We show that in the regular case, i.e., when 
all the  simplices~$\Delta_k$ are nondegenerate and 
the barycentric coordinates of~$\nill$ in $\Delta_k$ are bounded away from zero independently of~$k$, 
Algorithm~1 converges with a linear rate. Otherwise, if~$\Delta_k$ is close to degenerate in some iterations, we apply a  regularization procedure given in the next section (Algorithm~2). 
We begin Algorithm~1  with several auxiliary results. 
\bigskip
 
\subsection{Preliminary facts}

For an arbitrary vector~$\bb \in \re^n$, its {\em negative extension} 
is~$\{\, \lambda \, \bb\ : \ 
\lambda \le 0\ \}$. Thus, the negative extension of a nonzero vector is 
a ray directed opposite to~$\bb$; for a zero vector, it is zero.
\smallskip

A system of $N$ vectors $\bb_1, \ldots , \bb_{N} \in \re^n, \ N \ge n$,   
is said to be {\em nondegenerate} if every~$n$ of its vectors 
are linearly independent. Clearly, an independent system cannot contain zero vectors 
or a subset of~$m \le n-1$ dependent vectors. 
\begin{lemma}\label{l.10}
	Suppose~the convex hull  of vectors~$\bb_1, \ldots , \bb_{n+1} \in \re^n$ contains the origin; then for an arbitrary vector~$\bb_0 \in \re^n$, there exists~$s\in \{1, \ldots , n+1\}$, such that after the replacement of~$\bb_s$ by~$\bb_0$, the convex hull still contains the origin. 
	Moreover, if the enlarged system~$\{\bb_i\}_{i=0}^{n+1}$ is nondegenerate, then such an~$s$ is unique.  
\end{lemma}
{\tt Proof}. First consider the case when~$\{\bb_i\}_{i=0}^{n+1}$ is nondegenerate. 
Then  the 
simplex~$\Delta \, = \, {\rm co}\, \{\bb_s\}_{i=1}^{n+1}$
contains~$\nill$ in its interior. Hence,  the negative extension of~$\bb_0$ intersects a unique facet of~$\Delta$. Let~$\bb_s$ be the vertex   opposite to that facet;  then~$s$ is the desired index. Consider now a general system~$\{\bb_i\}_{i=0}^{n+1}$, which is admissible, i.e., 
satisfies the assumptions of the lemma. It can be slightly perturbed to become nondegenerate
and still be admissible. 
Hence, there exists a sequence of nondegenerate admissible 
systems~$\{\bb_i^{(k)}\}_{i=0}^{n+1}$
that converges pointwise to~$\{\bb_i\}_{i=0}^{n+1}$. By what we have proved above, 
for every~$k$, there exist a unique index~$s = s(k)$ such that the simplex with vertices 
$\{\bb_i^{(k)}\}_{i\ne s}$ contains the origin. Passing to a subsequence, it can be 
assumed that~$s(k)$ is the same for all~$k$. Now by passing to the limit we conclude that 
the simplex with vertices 
$\{\bb_i\}_{i\ne s}$ contains the origin. 
\bigskip

\subsection{The scheme of Algorithm~1}

Here we give a brief description of the algorithm and the
details will follow in subsection~4.4. We assume that there are 
points~$t_i\in K$ such that the system of moment vectors~$\{\bu(t_i)\}_{i=1}^{n+1}$
is nondegenerate. This is a part of the regularity assumption formulated in the next subsection. 
\bigskip 

\noindent {\em Input}. A system~$\Phi = \{\varphi_1, \ldots , \varphi_n\}$, a function~$f\in C(K)$, and $\varepsilon > 0$.
\smallskip 

\noindent {\em Output}. Numbers~$b, B$ and a polynomial~$\hat p \in \cP$
such that~$B - b < \varepsilon$, $\ b \, \le \, {\rm dist}\, (f, \cP) \, \le \, B$ and 
$\, b\, \le \, \|\hat p - f\| \, \le \, B$.  
\medskip 

\noindent {\em Initialization}. 
We start with arbitrary points~$t_1, \ldots , t_{n+1}$ for which the set of vectors~$\{\bu(t_i)\}_{i=1}^{n+1}$ is nondegenerate 
and choose the signs~$\sigma_i$ so that~the convex 
hull of the set~$\{\, \sigma_i\bu(t_i)\, \}_{i=1}^{n+1}$ contains the origin. 
Denote~$\ba_i = \sigma_i\bu(t_i)$ and 
find~$p \, \in \, \cP$ and~$d\in \re$ such that
the vector~$(p^1, \ldots , p^n, d)^T,$ where~$\bp = (p^1, \ldots , p^n)^T$ are coefficients of~$p$,  is a unique solution 
of the linear system 
\begin{equation}\label{eq.linear-system}
	\bigl(\bp, \bu(t_i)\bigr) \ = \ f(t_i)\ + \  
	\sigma_i \, d\, . 
\end{equation}
If~$d<0$, then change all the signs~$\sigma_i$. 
Thus, $d\ge 0$.  Set~$b_1 = d, \, B_1 = \|p-f\|_{C(K)}$. 
\bigskip 

\begin{figure}[h]
	    \caption{}
	    \centering
		\includegraphics[width=\linewidth]{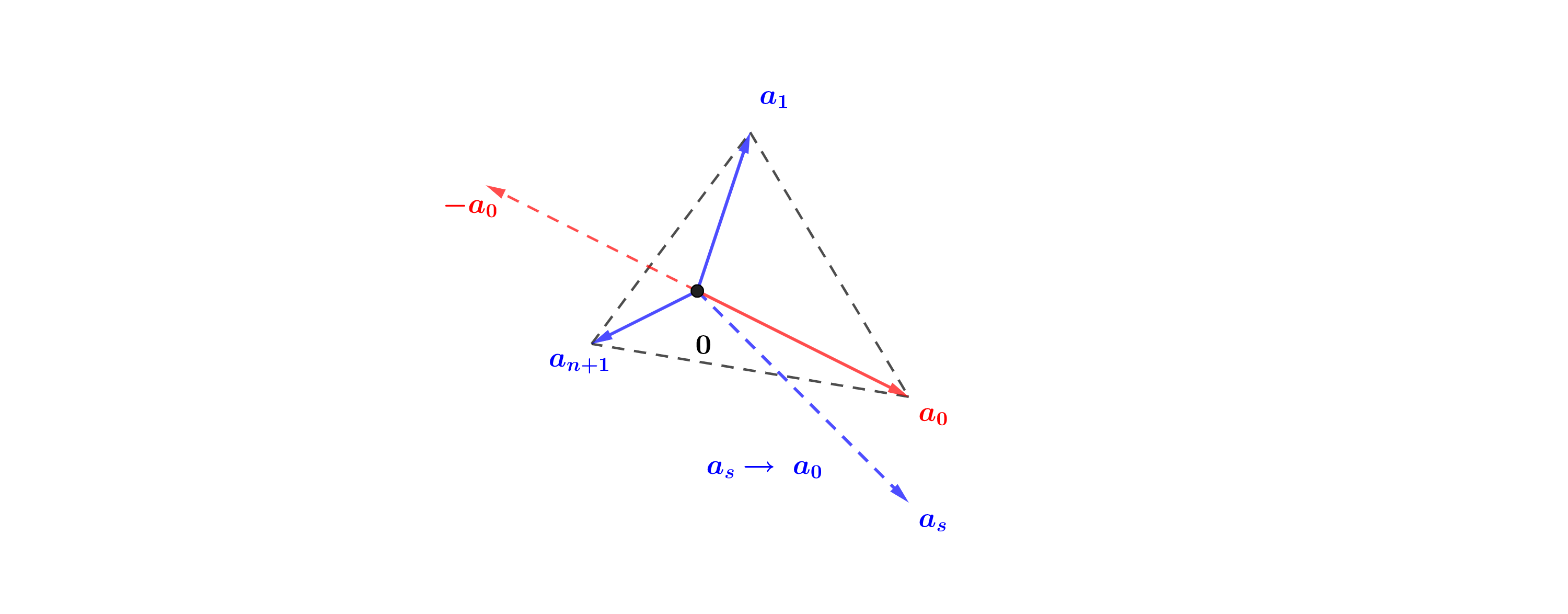}
	
	  \label{fig-algorithm}
\end{figure}

\noindent {\em Main loop}. After~$k-1$ iterations we have   
$n+1$ points (for the sake of simplicity, denote them as above~$t_1, \ldots , t_{n+1}$)
and the corresponding vectors~$\ba_1, \ldots , \ba_{n+1}$, whose convex hull contains the origin as an interior point. Then we take a point~$t_0$ of maximum of 
the function~$|p(t) - f(t)|$ on~$K$  and replace one of the vectors 
$\ba_s$ with~$\ba_0 = \sigma_0\bu(t_0)$ so that the new system still contains the origin inside its convex hull, see Fig.\ref{fig-algorithm}. 

 This means that the negative extension of~$\ba_0$
intersects the facet of the simplex~$\{\ba_i\}_{i=1}^{n+1}$
opposite to the vertex~$\ba_s$. If  the extended  system~$\{\ba_i\}_{i=0}^{n+1}$ is nondegenerate, then the index~$s$ is uniquely defined. Then we 
solve system~(\ref{eq.linear-system}) and 
set $b_{k+1} = \, d, \ B_{k+1} \, = \, \min\, \{B_k, \|p_k - f\|\}$.  
\bigskip 

\noindent {\em  Termination}. The algorithm halts when~$B_N - b_N < \varepsilon$.
Then we set~$b=b_N, \, B=B_N$ and~$\hat p = p_k$, where 
$k \le N$ is the number of iteration, for which~$\|p_k-f\| = B$. 
\bigskip

\subsection{The regularity assumption}

We say that the {\em regular case} takes place for Algorithm~1 if 
in each iteration, the system~$\{\ba_i\}^{n+1}_{i=0}$ is nondegenerate 
and the coefficient~$\alpha_0$ in the convex 
combination~$\nill \, = \, \sum_{i}\alpha_i \ba_i$ 
(the summation is over all~$i=0, \ldots , n+1$ except for~$s$) is not less than 
a prescribed number~$\mu > 0$. This number will be referred to as a {\em regularity parameter}.  
\bigskip

\subsection{The routine of Algorithm~1}

Here we give a detailed description of all steps. 
\medskip 

\noindent \textbf{Algorithm~1}. 
\smallskip 

\noindent {\em Input}. A system~$\Phi$ a function~$f\in C(K)$ and $\varepsilon > 0$.
\smallskip 

\noindent {\em Output}. Numbers~$b, B$ and a polynomial~$\hat p \in \cP$
such that~$B - b < \varepsilon$, $\ b \, \le \, {\rm dist}\, (f, \cP) \, \le \, B$ and 
$b\, \le \, \|\hat p - f\| \, \le \, B$.  
\medskip

\noindent  {\em Initialization.} 
Take an arbitrary set~$\cT_1\ = \, \{t_i\}_{i=1}^{n+1}$ of $n+1$ different points on~$K$ 
such that the system of moment vectors~$\{\bu(t_i)\}_{t_i\in \cT_1}$ is nondegenerate. 
Then the system of~$n+1$ linear equations~$\, \sum_{i=1}^{n+1} x_i \bu(t_i) \, = \, 0$
possesses  a unique solution~$\bx = (x_1, \ldots , x_{n+1})$. 
If~$\sigma_i$ is a sign of~$x_i$, $\alpha_i = |x_i|$ and $\ba_i = \sigma_i \bu(t_i)$, 
then~
$$
\sum_{i=1}^{n+1} x_i \bu(t_i) \ = \ \sum_{i=1}^{n+1} |x_i| \sigma_i \bu(t_i) \ = \ 
\sum_{i=1}^{n+1} |x_i| \ba_i \ = \ \sum_{i=1}^{n+1} \alpha_i\, \ba_i.
$$
Now find a vector~$\bp \in \re^{n+1}$
and a number~$d$ from the linear system~(\ref{eq.linear-system}). 
The corresponding polynomial~$p$ possesses the 
property~$p(t_i) \, = \, f(t_i) \, + \, \sigma_i d, \, i = 1, \ldots , n+1$.  
If~$d< 0$, then change all signs~$\sigma_i$. Thus, $b_1\ge 0$. Set~$p_1 = p, \, 
b_1 = d \, , \ B_1 = \|p_1-f\|$.  

\medskip 

{\em Main loop.} {\em The $k$th iteration}. 
We have a polynomial~$p_k \in \cP$,  a set of~$n+1$ points~$\cT_{k}\, = \, \{t_i^{(k)}\}_{i=1}^{n+1}$ (for the sake of simplicity we denote~$t_i^{(k)} = t_i$) 
and two numbers~$B_k > b_k \ge 0$
such that: 
\smallskip 

\noindent 1) $\ p_k(t_i)\, = \, f(t_k)\, + \, \sigma_i b_k$, where $\sigma_i \in \{-1, 1\}$ for all~$i$; 
\smallskip 

\noindent 2) 
the simplex~$\Delta_k \, = \, {\rm co}\, \{\ba_i\}_{i=1}^{n+1}$, where 
$\ba_i\, = \, \sigma_{i}\bu(t_i)$ contains the origin  
inside: $\nill \, = \, \sum_{i=1}^{n+1}\, \alpha_i^{(k)} \, \ba_i$ 
with all~$\alpha_i^{(k)}$ positive 
and~$\sum_i \alpha_i^{(k)} \, = \, 1$.   
\smallskip 

Denote~$r_k = \|p_{k} - f\|, \ \alpha_i \, = \, \alpha_i^{(k)}$, 
and find a point~$t_0$ for which~$|p_{k}(t_0) - f(t_0)| \, = \, r_{k}$. 
Clearly, $r_k \ge b_k$. Denote by~$\sigma_0$ the sign of~$p_{k}(t_0) - f(t_0)$ and
$\ba_0  = \sigma_0\bu(t_0)$. 

Now find the facet~$\Delta^s_k$ of~$\Delta_k$ intersected by the negative 
extension of~$\ba_0$. To this end, we solve the linear system~$\ba_0 \, = \, \sum_{i=0}^{n+1} x_i \ba_i$ 
and denote by~$\tau = \frac{x_s}{\alpha_s}$, the maximal ratio~$\frac{x_i}{\alpha_i}$. 
Then~$\ba_0 \, = \, \sum_{i=0}^{n+1} (x_i -  \tau\alpha_i) \ba_i$.
The coefficient~$x_s -  \tau\alpha_s$ vanishes and all other coefficients 
~$x_i -  \tau\alpha_i$ are negative. Hence, the negative extension of~$\ba_0$ intersects
the facet~$\Delta^s_k$ opposite to~$\ba_s$.  

Thus, replacing~$\ba_s$ by~$\ba_0$ we obtain a new simplex~$\Delta_{k+1}$ containing the origin. 
Replacing~$t_s$ by~$t_0$ we obtain a new system of points~$\cT_{k+1}$. Then we 
find 
a polynomial~$p_{k+1}$ and $b_{k+1} > 0$ from the linear system~$p_{k+1} (t_i) \, = \, f(t_i) +\sigma_i b_{k+1}, \ i \, =\, 0, \ldots , n+1, \, i\ne s$.
Set~$B_{k+1}\, = \, \min\, \{r_{k}, B_k\}$. 

Finally, redenote the points of the set~$\cT_{k+1}$ by $\{t_i\}_{i=1}^{n+1}$
in an arbitrary order. 

\medskip 

\noindent  {\em Termination.} If~$B_N - b_N < \varepsilon$, then STOP. 
\smallskip 

\noindent Set~$b=b_N, \, B=B_N$ and~$\hat p = p_k$, where 
$k \le N$ is the number of iteration, for which~${\|p_k-f\| = B}$. 
\bigskip

\noindent END

\medskip 

\begin{remark}\label{r.33}
	{\em The procedure of the Algorithm~1 stays  well-defined also without the 
		regularity assumption. Lemma~\ref{l.10} provides the existence 
		of the required point~$\ba_s$ even if the simplex~$\Delta_k$ is degenerate. 
		However, in this case the difference~$B_k - b_k$ may 
		not decrease and the termination criterion $B_N - b_N < \varepsilon$
		may not be reached within finite time. 
		
	}
\end{remark}

\bigskip

\subsection{The rate of convergence}

In the regular case Algorithm~1 converges with a linear rate, which 
follows from Theorem~\ref{th.30} below. 
\begin{theorem}\label{th.30}
	In the $k$th iteration of Algorithm~1, we have 
	\begin{equation}\label{eq.estim}
		B_{k+1}-b_{k+1} \ \le  \ \Bigl( 1 \, - \, \alpha_{0}^{(k)}\Bigr)\, 
		\Bigl(B_{k} \ - \ b_{k}\Bigr)\, , 
	\end{equation}
	where~$\alpha_0^{(k)} > 0$ is the coefficient in the 
	decomposition~$\nill \, = \, \sum\limits_{i}\, \alpha_i^{(k)} \, \ba_i$ 
	in the $k$th  iteration (the summation is over all~$i=0, \ldots , n+1$
	except for~$s$). 
\end{theorem}
{\tt Proof}. Without loss of generality it can be assumed that $s=n+1$, 
i.e., the~$k$th iteration replaces~$\ba_{n+1}$ with~$\ba_0$. Denote~$\alpha_i \, = \, 
\alpha_i^{(k)}$. 
We have $(\bp_k, \ba_i)\, = \, |p_k(t_i)|\, = \, b_k$
for all~$i=1, \ldots , n$ and 
$(\bp_k, \ba_0)\, = \, |\bp_k(t_0)|\, = \, r_k$. 
On the other hand~$\ (\bp_{k+1}, \ba_i)\, = \, |p_{k+1}(t_i)|\, = \, b_{k+1}, \ 
i = 0, \ldots, n$. 
Taking the difference, we obtain~$(\bp_{k+1} - \bp_{k}, \ba_i)\, = \, 
b_{k+1}-b_{k}$ for all~$i=1, \ldots , n$ and 
$(\bp_{k+1} - \bp_{k}, \ba_0)\, = \,  b_{k+1} - r_k$. 
Multiplying this equality by~$\alpha_i$, computing the sum over all~$i\le n$, 
we get 
\begin{multline*}
0 \ = \ 
\sum_{i = 0}^n \, \alpha_i \, \bigl(\bp_{k} - \bp_{k+1}\, , \, \ba_i\bigr) \ =\ 
\alpha_0(r_k - b_{k+1}) \ + \ \sum_{i = 1}^n \, \alpha_i (b_{k}-b_{k+1}) \ = \\ \ =
\alpha_0r_k \, - \, b_{k+1} \, + \ (1-\alpha_0) b_k\, . 
\end{multline*}
Thus, $\, \alpha_0r_k \, - \, b_{k+1} \, +  \,  (1-\alpha_0) b_k\, = \, 0$. 
Adding and subtracting $(1-\alpha_0)B_k$ we obtain after simple manipulations: 
\begin{equation}\label{eq.100}
	(1-\alpha_0)B_k  \, + \, \alpha_0 r_k \, - \, b_{k+1}\ = \ (1-\alpha_0) (B_k - b_k)
\end{equation}
Since~$B_{k+1} \, = \, \min \, \{B_k, r_k\}\ \le \ 
(1-\alpha_0)B_k  \, + \, \alpha_0 r_k$, we see that the left hand side of~(\ref{eq.100})
is bigger than or equal to~$B_{k+1} - b_{k+1}$, which concludes the proof.

{\hfill $\Box$}
\medskip 

Let us recall the regularity assumption: in all iterations we have~$\alpha_0 \ge \mu$, 
where $\mu > 0$ is fixed. 
\begin{theorem}\label{th.40}
	In the regular case Algorithm~1 converges linearly with the regularity parameter~$\mu$: 
	$$
	B_{k+1} \, - \, b_{k+1} \ \le \ (1-\mu)^k \, \bigl(B_{1} - b_{1}\bigr)\, , \quad k \in \n\, . 
	$$
\end{theorem}

\begin{remark}\label{r.35}
	{\em Equation~(\ref{eq.100}) shows that 
		if~$\alpha_0$ is small, then the $k$th iteration only slightly changes the difference~$B_{k} - b_k$. The algorithm slows down. 
		This happens if the vectors~$\ba_1, \ldots , \ba_n$ are almost linearly dependent and the linear systems~$\ba_0 \, = \, \sum_{i=0}^{n+1} x_i \ba_i$  becomes ill-conditioned.   Those problems can be resolved by the 
	 regularization procedure presented in the next section.}
\end{remark}

\begin{remark}\label{r.37}
	{\em From the proof of Theorem~\ref{th.30} we see that 
		it is not actually necessary to find the point of maximum~$t_0$
		at each iteration, it suffices to choose an arbitrary  point where the value~$|p_k(t_0) - f(t_0)|$  is significantly bigger than that at~$t_k$. This may be done, for instance, by a random search.  This strategy can be applied when finding the point of maximum of~$|p_k(t)-f(t)|$ is difficult.  
		For example, this situation occurs in the multivariate case. 
		
	}
\end{remark}
\bigskip

\subsection{Algorithm in the general case. Regularization}

The regularization procedure is needed when the simplex~$\Delta_k$
becomes close to degenerate. In this case we apply an essentially different 
method of choosing the next point in each iteration. This can be realized under additional assumptions:
\bigskip 

\noindent \textbf{Assumptions}. {\em $K \subset \re^N$ is a convex body,  
	all the functions 
	$\varphi_i$ are analytic on~$K$.}
\smallskip 

We begin with auxiliary results.  
In this section all lengths and distances in~$\re^d$ are measured in  the Euclidean norm, 
which is denoted by~$\|\cdot \|_2$.  

\begin{lemma}\label{l.20}
	Let vectors~$\bu(t_1), \ldots , \bu(t_{n-1})$ span a 
	hyperplane~$H \subset \re^n$. Then there exists a unique, up to the sign,  polynomial 
	$q \in \cP$ normalized by the condition~$\|\bq\|_2 = 1$
	which vanishes at the points~$t_1, \ldots , t_{n-1}$. 
	For  every~$t\in K$, the distance from~$\bu(t)$ to~$H$ is equal to~$|q(t)|$. 
\end{lemma}
{\tt Proof}. If~$\bq$ is a unit normal vector to~$H$, then
for the corresponding polynomial~$q$, we have~$q(t_i)\, = \, \bigl( \bq , \bu(t_i)\bigr) \, = \, 0$. If another polynomial~$y(\vardot)$ vanishes at the points~$t_i$, 
then its vector~$\by$ is orthogonal to all~$\bu(t_i)$ and hence, is 
orthogonal to~$H$, i.e., collinear to~$\bq$. Thus, the polynomial~$q$ is unique up to the sign.
The distance from~$\bu(t)$ to $H$ is equal to~$\, \bigl|\bigl( \bq, \bu(t)\bigr)\bigr| \, = \, 
|q(t)|$.

{\hfill $\Box$}
\medskip 

Lemma~\ref{l.20} expresses the distance from a moment vector to a hyperplane in terms 
of a special polynomial~$q$. This trick will be  applied in the regularization method. 

Let a set of  points~$\cT = \{t_i\}_{i=1}^{n+1}$ be such that the 
system of vectors~$\{\ba_i\}_{i=1}^{n+1}$ is nondegenerate, 
$p(t_i)-f(t_i)\,  = \, \sigma_i d$ for all~$i$, where~$d\ge 0$,   and $r= \|p-f\|_{C(K)}$.
For each pair~$i_1, i_2$ such that~$1 \le i_1 < i_2 \le n + 1$, 
let~$q_{i_1 i_2} \in \cP$ be the polynomial such that~$\|\bq_{i_1 i_2}\|_2 = 1$ 
and~$q_{i_1 i_2}(t_s) = 0$
for every~$s \notin \{i_1, i_2\}$. This polynomial is unique up to the sign.   
For a given~$\nu \in [0,1)$, we denote by~$t_{\nu}$ the solution of the following optimisation problem
\begin{equation}\label{eq.prob10}
	\left\{
	\begin{array}{l}
		g(t)\, = \, \sum\limits_{i_1 < i_2} \frac{1}{|q_{i_1 i_2}(t)|^2} \ \to \ \min\, , \\
		${}$\\ 
		\bigl|p(t) - f(t)\bigr|\, \ge \, (1-\nu)\, r \, + \, \nu \, d\, , \quad t\in K\,  \, .
	\end{array}
	\right. 
\end{equation}
\bigskip 
\begin{theorem}\label{th.50}
	For~$\nu \in [0, 1)$, let~$\rho(\nu)$ be 
	the maximal number such that there exists~$t\in K$ 
	for which $\bigl|p(t) - f(t)\bigr|\, \ge \, (1-\nu) r \, + \, \nu \, d$
	and the point~$\bu(t)$ is bounded away by the distance of at least~$\rho(\nu)$ from each of~$\frac{n(n+1)}{2}$
	hyperplanes spanned by~$n-1$ vectors 
	from~$\{\ba_{i}\}_{i=1}^{n+1}$ . 
	Then
	\begin{equation}\label{eq.g-rho}
		\frac{1}{\sqrt{g(t_{\nu})}} \quad < \quad  \rho(\nu) \quad \le \quad  \sqrt{\frac{n(n+1)}{2}}
		\frac{1}{\sqrt{g(t_{\nu})}}
	\end{equation}
	and the system~$\{\ba_{1}, \ldots , \ba_{n+1}, \bu(t_{\nu})\}$ is nondegenerate. 
\end{theorem}
{\tt Proof}. By Lemma~\ref{l.20} the distance from 
$\bu(t)$ to the hyperplane spanned by all vectors~$\bu(t_i)$
except for~$\bu(i_1), \bu(i_1)$ is equal to~$|q_{i_1i_2}(t)|$.
If we replace all~$\bu(t_i)$ by~$\ba(t_i)$, nothing changes since 
those vectors span the same hyperplane. 
If there exists an admissible~$t\in K$ for which $\bu(t)$ 
is bounded away from all those hyperplanes by the distance at least~$\rho$, then 
$g(t) \, \le \, \frac{n(n+1)}{2 \rho^2}$. Hence, $\rho \, \le \, 
\sqrt{\frac{n(n+1)}{2}}\frac{1}{\sqrt{g(t)}} \, \le \, \sqrt{\frac{n(n+1)}{2}}
\frac{1}{\sqrt{g(t_{\nu})}}$, which proves the upper bound in~(\ref{eq.g-rho}). 
Furthermore, $g(t_{\nu}) \, > \, \max\limits_{i_1 < i_2}\, \frac{1}{|q_{i_1 i_2}(t_{\nu})|^2}$, hence
the minimal number~$|q_{i_1 i_2}(t_{\nu})|$ over all~$i_1, i_2$
is greater  than~$
\frac{1}{\sqrt{g(t_{\nu})}}$. In view of Lemma~\ref{l.20}, 
this means that
the minimal distance from~$\bu(t_{\nu})$ to those hyperplanes is greater 
than~$\frac{1}{\sqrt{g(t_{\nu})}}$, which completes the proof of~(\ref{eq.g-rho}). 

If the system~$\{\ba_{1}, \ldots , \ba_{n+1}, \bu(t_{\nu})\}$ is degenerate, then it contains 
a subset of $n$ linearly  dependent vectors, one of which must be~$\bu(t)$, because 
the system~$\{\ba_i\}_{i=1}^{n+1}$ is nondegenerate. Hence, $\bu(t_{\nu})$ must belong to 
the linear span of other~$n-1$ vectors of this subset, in which case~$g(t_{\nu}) = + \infty$.

{\hfill $\Box$}
\bigskip 

Now we can present Algorithm~2 based on the aforementioned regularization. 
\smallskip  

\bigskip 

\noindent \textbf{Algorithm~2} coincides with Algorithm~1 with one crucial modification: 
we fix $\nu \in [0,1)$,  in each iteration solve problem~(\ref{eq.prob10}) 
for the set of points~$\cT = \cT_k$ and take  the optimal point~$t_{\nu}$
instead of~$t_0$ (the point of maximum of~$|p(t) - f(t)|$. Respectively, $r_k$ is replaced by~$(1-\nu)r_{k} + \nu b_k$ 

\smallskip 

The regular case is defined for Algorithm~2 in the same way as for~Algorithm~1: 
in all iterations~$\alpha_0 \ge \mu$. The rate of convergence in the regular case stays linear but slows down since 
the value~$r_k$ is replaced by a smaller value~$(1-\nu)r_{k} + \nu b_k$. Therefore, 
inequality~(\ref{eq.estim}) is replaced by
\begin{equation}\label{eq.estim-g}
	B_{k+1}-b_{k+1} \ \le  \ \Bigl( 1 \, - \, (1-\nu)\alpha_0\Bigr)\, 
	\Bigl(B_{k} \ - \ b_{k}\Bigr). 
\end{equation}
Consequently, the rate of linear convergence becomes 
$$
B_{k+1} \, - \, b_{k+1} \ \le \ \bigl(1\, -\, (1-\nu)\mu \bigr)^k \, \bigl(B_{1} - b_{1}\bigr)\, , \quad k \in \n\, . 
$$
In the regular case one should apply Algorithm~1, which is faster and less 
expensive.  If it slows down 
at degenerate simplices on several successive iterations, one should switch 
to Algorithm~2. 
\begin{remark}\label{r.50}
	{\em The key difference between Algorithms~1 and 2 is the choice of the new point 
		on each iteration. Algorithm~1 takes the point~$t_0$ of the maximal distance~$r_k = |p_k(t_0) - f(t_0)|$. Algorithm~2 instead considers all points~$t$, where this distance is  big enough, i.e., 
		$\, |p_k(t) - f(t)| \, \ge \, (1-\nu)r_k + d_k\, $. 
		On this set we find the point of minimum~$t_{\nu}$ of the 
		function~$g(t)$.  In other words, Algorithm~2 aims to guarantee the nongeneracy of the newly constructed simplex~$\Delta_{k+1}$ provided 
		that the value~$|p_k(t_{\nu}) - f(t_{\nu})|$  was significaltly larger than~$b_k$. The nondegeneracy, in turn, will guarantee 
		a sufficient descent of the value~$B_k - b_k$ on the next iterations. }
\end{remark}
\smallskip 

The value~$g(t_{\nu})$ can be estimated from above. Under our assumptions (all 
$\varphi_i$ are analytic and independent, $K$ is a convex body, $f$ is Lipschitz), 
it has the power growth in~$1/|r_k - d_k|$. 
\begin{prop}\label{p.20}
	We have  $\, g(t_{\nu}) \, = \, 
	O(|r_k - d_k|^{-2m})$, where~$m$ is the maximal multiplicity of zeros of polynomials from~$\cP$
\end{prop} 
{\tt Proof}. Suppose~$M = \max_{p\in \cP, \|p - f\| \le r_k}\|p'\|_{C(K)}$ and $L$ is the Lipschitz constant of~$f$, then the inequality~$|p(t) - f(t)| \, \ge \,  \frac{1}{2} (r_k + d_k)$ is valid for all~$t$ from the ball~$B(t_0, R)$ of radius~$R=\frac{1}{2(M+L)} (r_k - d_k)$
centered at~$t_0$. On the other hand,  there are constants~$C, R_0 > 0$ such $|q_{i_1 i_2}(t)| \, > \, 
C\|t-t_0\|^{m}$, whenever~$\|t-t_0\| < R_0$. Then for every~$t \in K$
such that~$\|t-t_0\| < \min\{ R, R_0\}$, we have 
$g(t) \, = \, O(|r_k - d_k|^{-2m})$, which is required. 

{\hfill $\Box$}
\bigskip 

\smallskip 

\medskip 

\textbf{Practical issue}. In all numerical results we set~$\nu = 0.5$. Thus, 
$t_{0.5}$ is the point of minimum of~$g(t)$ 
on the set~$\, \bigl\{ t \in K \ : \ |p_k(t) - f(t)| \ge \frac{1}{2}(r_k + d_k) \, \bigr\}$.  
In practice it often suffices to apply a simple procedure, without solving any optimisation problem. We fix a small~$\delta > 0$ and at every iteration 
compute $n$ values:~$|q_{s, j}(t_0)|, \, j\ne s$, where~$t_0$ is the point of 
maximum of~$|p_k(t) - f(t)|$ on~$K$
and~$s$ is the index of the point~$t_s$ replaced by~$t_0$ in the $k$th iteration, 
$j = 1, \ldots , n+1$,  $j\ne s$. 
Let us recall that~$q_{s, j}$ is the polynomial with zeros at the points~$t_i, \, i\notin \{s, j\}$. If $|q_{s, j}(t_0)| \, < \, \delta$ for some~$j$, we find a point~$\bar t$
for which~$|p_k(\bar t) - f(\bar t)| \ge \frac12\bigl( r_k + d_k\bigr)$ 
and~$|q_{s, j}(t_0)| \, > \, \delta$. This point can be found, for example, 
by the penalty function method. 

\medskip 

\begin{remark}\label{r.60}
	{\em The regularization moves the point~$\ba_0$ away from all 
		hyperplanes spanned by the other $(n+1)$ vectors~$\ba_i$. 
		Hence, it always keeps the simplices~$\Delta_k$ nondegenerate. This, however, does not 
		imply that they are degenerate uniformly for all~$k$ and that the parameter~$\alpha_0^{(k)}$
		can be bounded away from zero by some positive constants. Although, a vast majority of the 
		numerical experiments demonstrate the uniform nondegeneracy of those simpleces, and, respectively, linear convergence of Algorithm~2. Nevertheless, we do not have a proof of its convergence for any initial data and can only formulate the corresponding conjecture.   
		According to numerical experiments, Algorithm~1 handles from~$20\%$ to~$80\%$
		(depending on the problem)
		cases without suffering by degeneracy. The remaining cases need the regularization
		and are always handled by  Algorithm~2.  We report the numerical results in Section~9.

	}
\end{remark}

\bigskip 

\begin{center}
	\large{\textbf{6. Approximation under linear constraints}}
\end{center}

\medskip 

Unlike the classical Remez algorithm, the method introduced in Sections 3 and 4 is easily extended  to constrained  approximation.   For   given linear functionals~$\ell_j$ and numbers~$b_j$, we 
consider the problem~(\ref{eq.constr}).  
As before, $f$ is continuous on a compact set~$K$, 
$\cP = {\rm span}\, \{\varphi_1, \ldots , \varphi_n\}$. 
Clearly, it can be assumed that all~$\ell_j$ are independent and therefore, $r \le n-1$. 

A Chebyshev system may loose its properties after imposing a linear constraint
(Example~\ref{ex.15}). However, the concept of generalized alternance works for the constrained problems as well.   

We denote~$L = \{\ell_j\}_{j=1}^r, \ 
L^{\perp} \, = \, \{p \in \cP\, : \, \ell_1(p) = \ldots = \ell_r(p) = 0\}$ and as usual identify~$\cP$ with~$\re^n$. Every functional~$\ell_j$ is associated to the vector~$\bell_j \in \re^n$
such that~$\ell_j(p) = (\bell_j, \bp)$. The orthogonal complement to the vectors~$\bell_1, \ldots , \bell_r$ will be denoted by the same  symbol~$L^{\perp}$. Since~$\ell_j$ are 
linearly   independent, it follows that~${\rm dim}\, L^{\perp} \, = \, n-r$. By projection we always mean an orthogonal projection. The  projections of~$\varphi_i$ to~$L^{\perp}$ are assumed to be independent, otherwise, we 
remove redundant functions from~$\Phi$. 

\begin{defi}\label{d.10c}
	Let~$\cP$ be the space of polynomials by a system~$\Phi = \{\varphi\}_{i=1}^n$, 
	$L = \{\ell_j\}_{j=1}^r, \, r< n$,  be a set of independent linear  functionals on~$\cP$. 
	For a given polynomial~$p\in \cP$,  a set of points~$\{\tau_1, \ldots , \tau_{m}\}\subset K, \ m\le n-r+ 1$, is called {\em $L$-alternance} if 
	\medskip 
	
	\noindent 1) for each~$i =1, \ldots , m$, we have $\, p(\tau_i) - f(\tau_i)\, = \, \sigma_i \, \|p-f\|, \, i=1, \ldots , m$, where $\sigma_i \in \{-1, 1\}$. 
	\smallskip 
	
	\noindent  2) the convex hull of the projections of the 
	vectors~$\{\sigma_i\bu(\tau_i)\}_{i=1}^m$ 
	onto~$L^{\perp}$ 
	contains the origin. 
\end{defi}
The proofs of the following generalization of Theorem~\ref{th.10} 
is literally the same:   
\begin{theorem}\label{th.10c}
	A polynomial~$p \in \cP$ is of the best approximation for a function~$f$ 
	under the constraints~$\ell_j(p) = b_j, \, j = 1, \ldots , r$, if and only if 
	it has an L-alternance. 
\end{theorem}
In case of non-uniqueness, all polynomials of the best approximation are characterised by 
a common~L-alternance. The proof of the following theorem is the same as for the 
corresponding Theorem~\ref{th.15}: 
\begin{theorem}\label{th.15c}
	For every~$f$, all polynomials of the best approximation from~$\cP$
	under constraints~$\ell_j(p) = b_j, \, \ell_j\in L$, have the same $L$-alternance 
	and the same values on it.  
\end{theorem}

Algorithms 1 and 2 are transferred to the  constrained problem~(\ref{eq.constr}) merely by replacing all vectors~$\ba_i$ by their projections to~$L^{\perp}$. Respectively, the dimension~$n$
is replaced by~$n-r$ and~$\re^n$ is replaced by~$L^{\perp}$. In Algorithm~2, the polynomial 
$q_{i_1 i_2}$ is defined by the set of zeros~$t_i, \ i \in \, \{1, \ldots , n-r+1\}\setminus 
\{i_1, i_2\}$
and by the conditions~$\ell_j(q_{i_1 i_2}) = 0, \, \ j=1, \ldots , r$. We refer to those algorithms as Algorithms~1c and 2c respectively (``c'' is from ``constrained''). 
\medskip 

\begin{remark}\label{r.70}
	{\em It is possible to obtain some generalizations of Theorem~\ref{th.10c} to nonlinear smooth 
		constraints~$\ell_j$. 
		In this case the subspace~$L$ will be spanned by the derivatives~$\ell_j'(p)$.
		The analogue of Theorem~\ref{th.10c} is valid in the direction ``necessity'': 
		every optimal polynomial do possess an $L$-alternance, but the converse is in general not true.   
		It becomes true if all~$\ell_j$ are convex and all equality constraints are replaced 
		by inequalities~$\ell_j(p) \le b_j$. In this case we have a complete analogue of  
		Theorem~\ref{th.10c} with the subspace~$L$ 
		spanned by the derivatives~$\ell_j'(p)$ of active constraints (computed at the point~$p$).

		However, for Algorithms 1c and 2c, the linearity of functionals~$\ell_j$ is significant. 
		The reason is that  
		if some of the functionals~$\ell_j(p)$ are nonlinear, then  the subspace~$L$  depends on~$p$. 
		Therefore, the subspaces~$L(p_k)$
		and $L(p_{k+1})$ in Algorithms 1c can  be different.  The key requirement 
		of the algorithm that the 
		projections of the simplex~$\Delta_k$ to~$L^{\perp}$ contains the origin 
		may not be preserved in $k$th iteration.

	}
\end{remark}

\bigskip

\bigskip

\begin{center}
	\large{\textbf{7. Approximations on unbounded domains}}
\end{center}

\bigskip

All our results are valid 
for arbitrary closed sets~$K$ (and convex in Algorithms 2 and 2c), not necessarily bounded, 
provided that all~$\varphi_i$ and~$f$ tend to zero as~$t \to \infty$. Indeed, 
if the distance $r = \min_{p\in \cP}\|p-f\|_{C(K)}$ is positive (otherwise the statement is trivial), then it is positive on some ball~$B \subset K$. 
The set of polynomials~$p$ for which~$\|p-f\|_{C(B)} \le  2d$ is compact in~$\cP$. 
Hence, all polynomials~$p$ uniformly  tend to zero
on this set  as~$t\to \infty$. Thus, there exists a ball~$B_1 \subset B$ such that
the best approximation polynomial on~$B_1$ coincides with that on~$K$. 
This reduces our problem to the compact set~$B_1$. 
The criteria of the best approximation in terms of the generalized alternance
and L-alternance, as well as the algorithm are extended to unbounded domain. 
The only assumption is that all~$\varphi_i$ decay at infinity.

\bigskip

\begin{center}
	\large{\textbf{8. Applications}}
\end{center}

\bigskip

\begin{center}
	\textbf{8.1. Structural recovery of signals}
\end{center}
\medskip

The problem of decomposing a continuous time signal~$f(t)$ 
into the form \\ $f(t) = \sum_{i=1}^n c_i \varphi(t, z_i) $, where 
$\varphi(t, z_i)$ are functions from a certain class,  arises in a number of applications such as sparse array of arrival estimation \cite{Kn}, system theory for parametrized model reduction \cite{Io}, radar \cite{Om}, ultrasound \cite{Ronen}, etc. For the case when $ \varphi(t,z_i) = e^{2 \pi z_i t} $ and $ z_i, c_i $ are unknown complex and real parameters respectively, there are many methods. A popular approach is related to Prony's method \cite{Berent, Beylkin, Plonka}. If $ \varphi(t,z_i) = \varphi(t - z_i)  $, for example $ \varphi(t,z_i) = e^{\frac{-(t-z_i)^2}{\sigma^2}} $ or $ \varphi(t,z_i) = \frac{1}{1+\frac{(t-z_i)^2}{\sigma^2}} $, where $ \sigma > 0 $,  the matrix pencil method is efficient~\cite{Yingbo}. Most of those methods are unstable in the presence of deterministic or random noise. In this case the decomposition problem is replaced by that of approximation and sometimes requires different tools. Algorithms 1 and 2 find the best approximation of a signal in the case when all the decomposition functions~$\varphi_{t, z_i}$ are known.

Consider a continuous time attenuating signal~$f\in C^1(\re_+)$  on the positive half-line, 
which tends to zero as~$t\to +\infty$. One needs to find its best approximation 
by a function of the form~$p(t)\, = \, \sum_{m}\, C_m(t)\, e^{\, z_mt}$
where~$z_m = a_m + b_mi, \ b_m<0$,  are known complex frequencies and~$C_m(t)$
are unknown algebraic polynomials of prescribed powers (they correspond to multiple frequencies).  We specialize our investigation to real-valued signals, when $p$ can be rewritten as a quasipolynomial~$p(t)\, = \, \sum_{k}\, P_k(t)\, e^{\, a_kt}\, \cos \, b_kt\, + \, Q_k(t)\, e^{\, a_kt}\, \sin \, b_kt$. Here $P_k, Q_k$ are real algebraic polynomials and 
${\sum_{k}\, {\rm deg}\, P_k(t) \, + \, {\rm deg}\, Q_k(t)\, = \, n}$. 
In most of applications all the polynomials are constants and one has to realize
an optimal  recovery 
in the uniform metric in~$C(\re_+)$ by finding coefficients~$c_m$. 

The optimal recovery problem is also considered with a sampling. 
For example, we know the values of~$p(t)$ at some points or other linear parameters (derivatives, Fourier coefficients, etc.) 
In this case we have a problem with linear constrains which can be solved 
with Algorithms 1c and 2c.  Numerical examples are given in Section~9. 
\smallskip

\bigskip

\begin{center}
	\textbf{8.2. Markov-Bernstein type inequalities}
\end{center}

\medskip

The Markov-Bernstein constants are  sharp upper bounds~$C_j = C_j(\Phi)$
in inequalities between the functions and their 
derivatives:~$\|p^{(j)}\|_{C(K)}\, \le \, C_j \, \|p\|_{C(K)}, \, p\in \cP$.
The classical Bernstein theorem asserts that for trigonometric polynomials 
$p = \sum_{k=0}^N\, (a_k \cos kt \, + \, b_k \sin kt)$ on the period~$[-\pi, \pi]$, 
the sharp constant~$C_1 = N$
(respectively, $C_j = N^j$) is attained on the functions~$p(t) = \cos Nt$
and $p(t) = \sin \, Nt$. The Markov inequality establishes that for 
algebraic polynomials of degree~$N$ on~$K = [-1, 1]$, the sharp constant is~$C_1 = N^{2}$ 
and 
it is   attained on the Chebyshev polynomial~$T_N \, = \, \cos \, (N\, {\rm acos}\, t)$. 
See~\cite{G} for a survey and extensions to entire functions. 
Many works are concerned with  the generalization of the Markov-Bernstein type inequalities 
for polynomials over other functional systems. In~\cite{Sz, MN} this was done 
for the system~$\{t^{\,k-1}e^{-t}\}_{k=1}^{n}$
(algebraic polynomials with the Laguerre weight), in~\cite{Fre} it is done for 
algebraic polynomials with a Gaussian weight. Systems of real negative exponents~$\{e^{\,a_kt}\}_{k=1}^{n}, \, a_k < 0, \, k = 1,\ldots , n$,  
have been analysed in~\cite{BE1, BE2, New, PJ1, P22}. 
Both these systems are Chebyshev. Generalizations of Markov-Bernstein inequalities  to non-Chebyshev systems, apart from the theoretical interest, is important 
in analytic geometry, ODEs, and dynamical systems~\cite{Fra, PJ1}. 
The latter application deals with the systems of complex exponentials in~$\re_+$~\cite{P22}.  

We begin with lacunary 
polynomials spanned by an 
arbitrary collection of integer 
nonnegative powers~$\varphi_k = t^{m_k}, \, k=1, \ldots , n$. They are usually 
non-Chebyshev.  
\begin{prop}\label{p.30}
	For lacunary algebraic polynomials \\ $\cP \, = \, \bigl\{p(t)  \, = \, 
	\sum_{k=1}^{n} p_k t^{\,m_k}\ : \ 
	p_k \in \re \bigr\}$ on~$K=[-1,1]$, the constant~$C_j$ is the reciprocal to the objective value of the problem
	\begin{equation}\label{eq.MB}
		\left\{
		\begin{array}{l}
			\|p\|_{C[-1,1]} \ \to \ \min\, \\
			${}$\\ 
			\ p^{(j)}(-1) \ = \ 1\, .
		\end{array}
		\right. 
	\end{equation}
\end{prop}
Thus, the Markov-Bernstein constants   for lacunary polynomials 
are obtained using problem~(\ref{eq.constr}) with~$f \equiv 0, \, \ell(p) = p^{(j)}(-1)$. 
\smallskip 

\noindent {\tt Proof.} The minimal norm~$\|p\|$ under the assumption~$\|p^{(j)}\| = 1$
is equal to~$1/C_j$. Let it be attained on some polynomial~$p$ and 
$\|p^{(j)}\| =  |p^{(j)}(a)|
$ at some point~$a\in [-1,1]$, say, $a < 0$. 
If~$a > - 1$, then the polynomial~$\bar p(t) \, = \, a^j  p(t/a)$ 
satisfies~$\|\bar p\| \, \le \, |a|^j  \|p\|$, while $\|\bar p^{(j)}\|\, = \, |\bar p^{(j)}(-1)| \, = \, |a|^j |a^{-j}p^{(j)}(a)|\, = \, 1$, which contradicts to the optimality of~$p$. 
Thus, $|p^{(d)}(-1)| = 1$, and hence, up to the sign, $p$ coincides with the 
solution of the aforementioned extremal problem.

{\hfill $\Box$}
\bigskip 

The second example is given by  real-valued polynomials of complex exponentials
(quasipolynomials). We assume that all 
frequencies~$z_k$ are of multiplicity one.  
\begin{prop}\label{p.40}
	For systems~$\varphi_k \, = \, e^{\,z_k t}, \, {\rm Re}\, z_k \, < \, 0, 
	\,$ $ k=1, \ldots , n$, on the half-line~$\re_+$, the constant~$C_j$ is the reciprocal to the objective value of the problem~$\|p\|_{C[0, +\infty]} \to \min\, , \ p^{(j)}(0) = 1$.  
\end{prop}
{\tt Proof.} If the minimal norm~$\|p\|$ under the assumption~$\|p^{(j)}\| = 1$
is attained on some polynomial~$p$, for which 
$\|p^{(j)}\| =  |p^{(j)}(a)|$, then 
for the polynomial~$\bar p(\vardot)\, = \, p(\vardot + a)$, we have 
~$\|\bar p\| \le \|p\|$ and~$\|p^{(j)}\|\, = \, |p^{(j)}(0)|$. Hence, up to the sign, $p$ coincides with the 
solution of that extremal problem.

{\hfill $\Box$}
\bigskip

\bigskip 

\begin{center}
	\textbf{8.3. Bounded trajectories of linear ODEs}
\end{center}

\medskip 

We consider a linear~ODE of $n$th order: $x^{(n)}\, = \, \sum_{k=0}^{n-1}c_kx^{(k)}$ 
on~$[0, +\infty)$ with 
constant coefficients~$c_k$ under the asymptotic stability assumption (every solution tends to zero as~$t\to \infty$). The problem is how large the derivative~$\dot x (0)$
can be if~$|x(t)| \le 1$ for all~$t\ge 0$.  
Each solution~$x(\cdot)$ is a linear combination of~$n$ elementary solutions~$t^se^{\, zt}$, where 
$z$ is a root of the characteristic polynomial~$P(x) \, = \, z^n \, - \, \sum_{k=0}^{n-1}c_k z^k$ 
and~$s$ does not exceed its multiplicity. Hence, the problem is equivalent to 
minimizing~$\|p\|_{C(\re_+)}$ under the constraint~$\dot x(0) = 1$. It it solved as 
in subsection~8.2 (computing~Markov-Bernstein constants).

\bigskip 

\begin{center}
	\textbf{8.4. Linear switching dynamical systems}
\end{center}
\medskip 

Linear switching system is a linear ODE  $\dot \bx (t) \, = \,  A(t)\bx(t)$
on the vector-function $\bx \, : \, \re_+ \, \to \, \re^n$ with the initial 
condition~$\bx(0) = \bx_0$ such that  the matrix function $A(t)$ 
takes values from a given compact set~$\cA$ called 
a {\em control set}. The 
{\em switching law}  is an piecewise-constant function  
$A: \re_+ \to \cA$. 
The {\em switching interval} is the interval between two 
consecutive switching points. This is the time of a continuous action of one matrix. 
The  {\em dwell time assumption} means that each interval is not less than a fixed time~$m$.

Linear switching systems have found many applications, mostly in electronic engineering, 
robotics,  and planning~\cite{B, GC, L}. 
One of the main issues  is the stability of the system~\cite{BCM, CGPS, MP}. 
A system is called {\em asymptotically stable}  if all its trajectories tend to zero. 
Deciding stability is usually hard~\cite{GO}. One of the problems in this direction is to 
establish the minimal stability interval~\cite{KP}. This is the smallest number~$M$ 
with the following property:  
if all trajectories with lengths of switching intervals at most~$M$ tend to zero, then the system is stable. In other words, it is always sufficient to check the stability 
only for switching laws with short intervals. It is known~\cite{PK} that~$M = m+T$, where
$T$ is the minimal number such that for every~$A \in \cA$, the value~$p(T)$ of the problem: 
\begin{equation}\label{eq.extr-cut}
	\left\{
	\begin{array}{l}
		\|\bp\|_{C(\re_+)} \to \min, \\ 
		p(T) = 1,\\
		p \in {\cP}_A
	\end{array}
	\right. 
\end{equation}
is bigger  than~$1$.  Here~$\cP_{A}$ denotes the space of polynomials 
defined by the system~$\Phi \, = \, \bigl\{\, t^ke^{\, zt}\ : \ z \in {\rm sp}\,(A), \ k \le 
d(z)\bigr\}$, where 
$ {\rm sp}\, (A)$ is the spectrum of~$A$ and~$d(z)$ is the size of the largest Jordan block 
of the eigenvalue~$z$. 

Thus, to find~$M$ one needs to solve problem~(\ref{eq.extr-cut}) for each~$A \in \cA$. 
Then by bisection we find the minimal number~$T$ for which the value of all problems is bigger 
than one.  In turn, problem~(\ref{eq.extr-cut}) is solved as the approximation 
problem~(\ref{eq.constr}) with one constraint~$\ell(p) = p(T) = 1$.  

\newpage

\begin{center}
	\large{\textbf{9. Examples and numerical results}}
\end{center}
\bigskip 

The code can be found on the website \smallskip 
\smallskip 
 
{https://gitlab.com/rinatik7.kamalov.7/approximation-by-non-chebyshev-systems-of-functions}
\smallskip 

{https://gitlab.com/users/rinatik7.kamalov.7/projects}.

\bigskip
\begin{center}
	\textbf{9.1. Approximation by shifts of Gaussian functions}
\end{center}
\bigskip

\begin{figure}[h]
	\caption{}
	\centering
	\begin{subfigure}{0.3\textwidth}
		\centering
		\includegraphics[width=\linewidth]{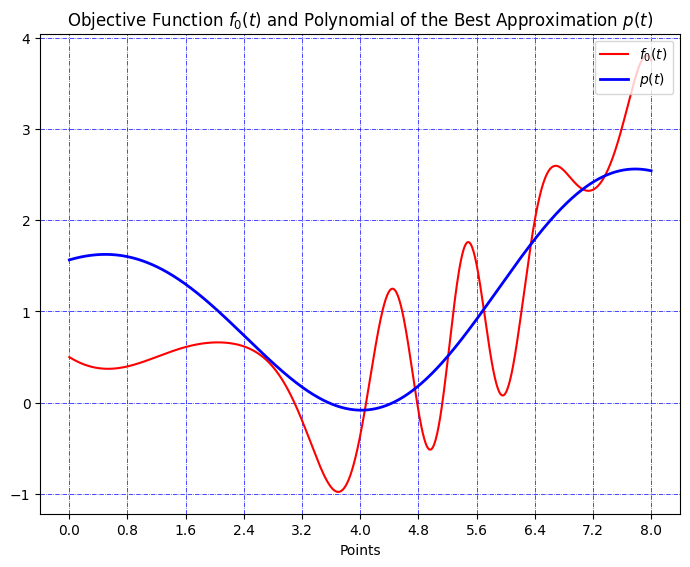}
		\caption{$ f(t) $ and $ p(t) $}
	\end{subfigure}\hfill
	\begin{subfigure}{0.3\textwidth}
		\centering
		\includegraphics[width=\linewidth]{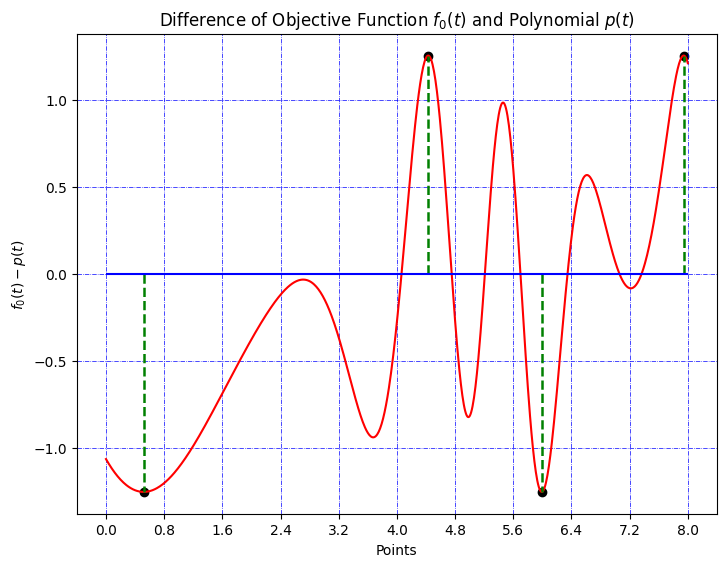}
		\caption{$ f(t) - p(t) $}
	\end{subfigure}\hfill
	\begin{subfigure}{0.3\textwidth}
		\centering
		\includegraphics[width=\linewidth]{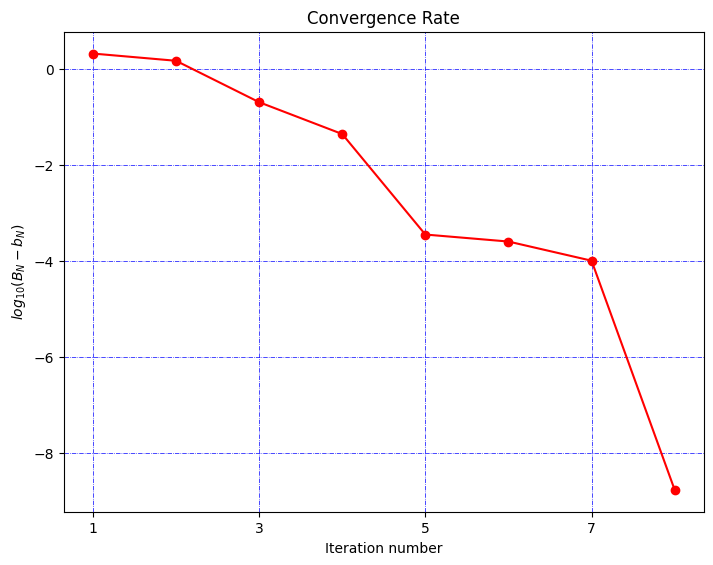}
		\caption{$ \log_{10} (B_N - b_N) $}
	\end{subfigure}
	\label{Fig:Gaussian1}
\end{figure}

\textbf{Unconstrained problems.} We approximate the signal 
\begin{equation*}
	f(t) \ = \ \frac{1}{10} \left(t-5 \right)^2 + \frac{1}{2} \left(t-4 \right) + \sin \left(0.4 \cdot t^2\cos \left(0.5 t\right )\right)\, , 
\end{equation*}
with polynomials by the 
system~$\Phi \, = \, \bigl\{   e^{\frac{-(t-1)^2}{9}}\, , \, e^{\frac{-(t-5)^2}{9}}\, , \,  e^{\frac{-(t-7)^2}{9}} \bigr\}$ on the segment $\, K = [0, 8] $ with $ \varepsilon = 10 ^ {-6} $. Eight iterations of the Algorithm~2 yield the required accuracy (Fig.\ref{Fig:Gaussian1}). The best approximation polynomial $ p(t) $ has corresponding coefficients  $\bp \, = \, (1.902091, -2.453699, 3.842463)^T$ with the 
distance  $  \ \| f - p \| \, = \,  1.254985 $. 
The generalized alternance consists of points
\, $\, \tau_1 = 0.517919, \, \tau_2 = 4.430493, \, \tau_3  =  5.992115, 
\tau_4  =  7.942944.$ 

\begin{figure}[h]
	\centering
	\caption{}
	\begin{subfigure}{0.3\textwidth}
		\centering
		\includegraphics[width=\linewidth]{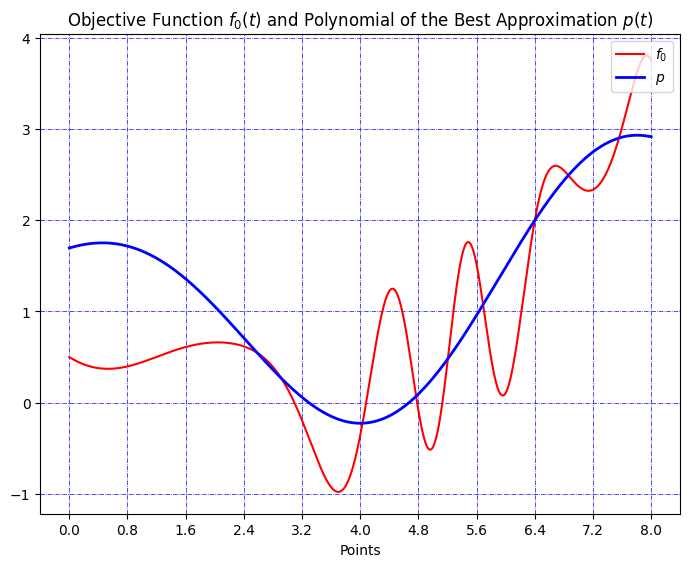}
		\caption{$ f(t) $ and $ p(t) $; $ p(6.4) = 2 $}
	\end{subfigure}\hfill
	\begin{subfigure}{0.3\textwidth}
		\centering
		\includegraphics[width=\linewidth]{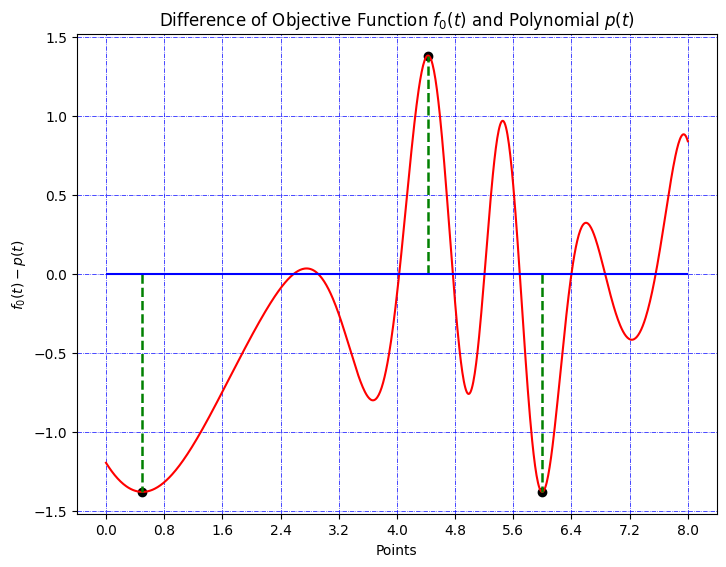}
		\caption{$ f(t) - p(t) $}
	\end{subfigure}\hfill
	\begin{subfigure}{0.3\textwidth}
		\centering
		\includegraphics[width=\linewidth]{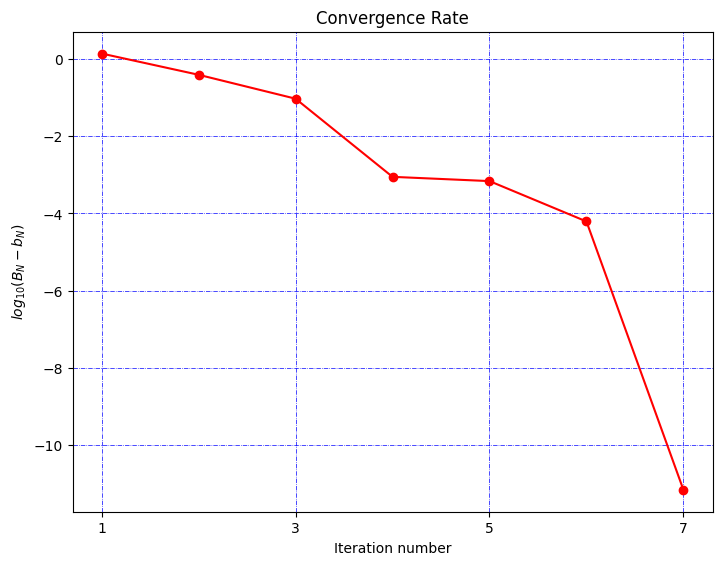}
		\caption{$ \log_{10} (B_N - b_N) $}
	\end{subfigure}
	\label{Fig:Gaussian2}
\end{figure}

\textbf{A problem with linear constraints.}  
In the previous example, the value of the approximated function $ f(t) $ 
at the point $t= 6.4 $ is  $ 1.999... $ while 
for the best approximation polynomial,  is~$p(6.4) = 1.78...$ If we need 
to approximate~$f$ and keep its value at~$t=6.4$, which we for simplicity round to~$2$,  we come to the 
constrained problem 
\begin{equation}
	\begin{cases}
		\label{approximation_problem_linear}
		\| f - p \|\, \to \, \min \\
		l(p) \, = \,  p(6.4)\, = \, 2\, \\
		p\in \cP
	\end{cases}
\end{equation}
The results of  Algorithm~2c are shown in the Fig.~\ref{Fig:Gaussian2}. 
The best approximation polynomial $ p(t) $ has corresponding coefficients  
$p_1 = 2.078450, \, p_2 \, = \, -2.939696, \, p_3\, = \, 4.457802$ with the 
minimal distance  $\| f - p \| \, = \,  1.3807$. The $L$-alternance comprises 
$\, \tau_1 = 0.500162, \, \tau_2 =  4.427931, \, \tau_3 = 5.998317.$
\medskip

Now add an extra condition:  the value of the derivative of $ p(t) $ at the same point 
$t=6.4$  coincides with that of~$f$, which is approximately  $ 4.47 $.
Thus, we add the constraint~$l_2(p) \, = \,  p'(6.4) = 4.47$. 

\begin{figure}[h]
	\centering
	\caption{}
	\adjustbox{valign=m}{\begin{subfigure}{0.3\textwidth}
			\centering
			\includegraphics[width=\linewidth]{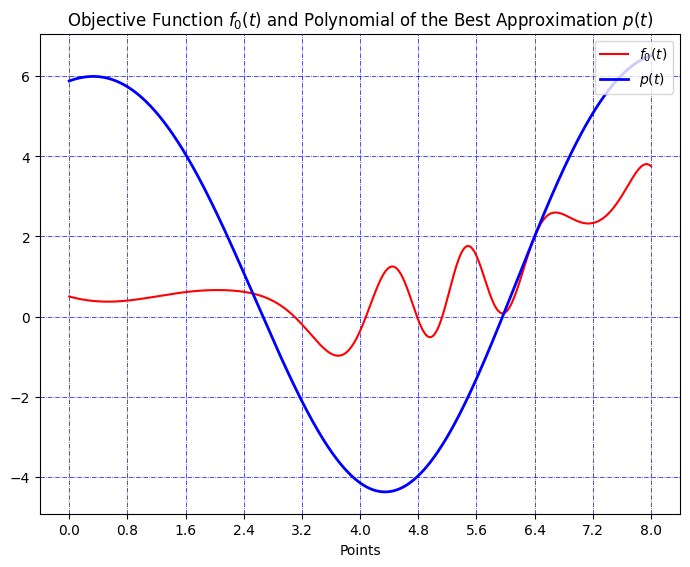}
			\caption{$ f(t) $ and $ p(t) $; \\ $ p(6.4) = 2, \,  p'(6.4) = 4.47 $}
	\end{subfigure}}\hfill
	\adjustbox{valign=m}{\begin{subfigure}{0.3\textwidth}
			\centering
			\includegraphics[width=\linewidth]{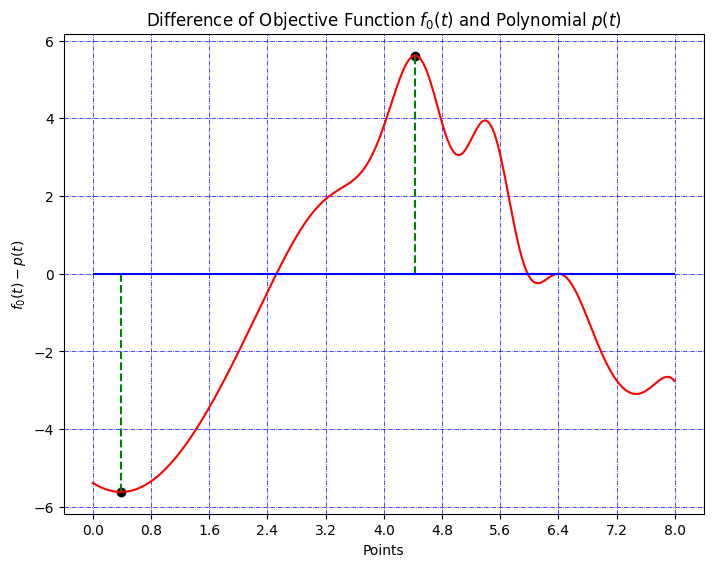}
			\caption{$ f(t) - p(t) $}
	\end{subfigure}}\hfill
	\adjustbox{valign=m}{\begin{subfigure}{0.3\textwidth}
			\centering
			\includegraphics[width=\linewidth]{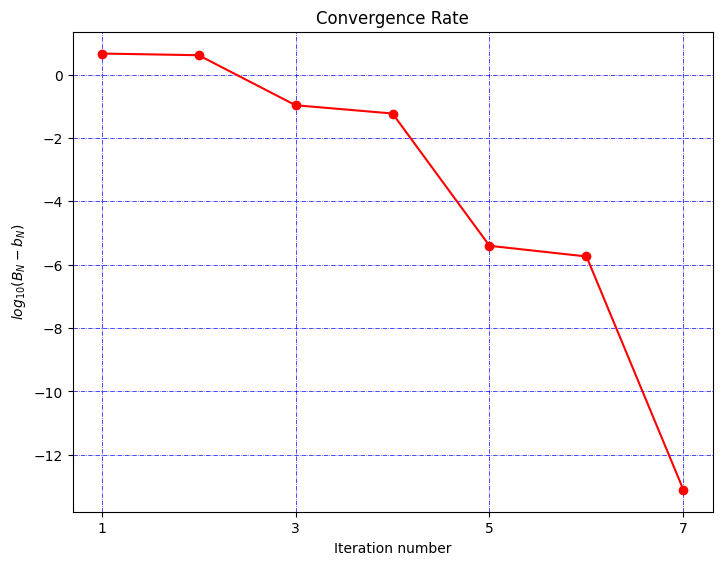}
			\caption{$ \log_{10} (B_N - b_N) $}
	\end{subfigure}}
	\label{Fig:Gaussian3}
\end{figure}

The Algorithm $2c$ produces  the following results (see Fig.\ref{Fig:Gaussian3}). The coefficients of the best polynomial  are   
$p_1 \, =  \, 7.407235, \, p_2 \, = \,  -12.84065, \, p_3 = 12.52896$ with the minimal distance $  \| f - p \|\, = \, 5.614225 $. The $L$-alternance consists of two points 
$\tau_1 = 4.430836, \, \tau_2 = 0.386453.$
\bigskip 

\begin{center}
	\textbf{9.2. Signal recovery by complex exponentials}
\end{center}
\bigskip

\begin{figure}[h]
	\caption{}
	\centering
	\adjustbox{valign=m}{\begin{subfigure}{0.3\textwidth}
			\centering
			\includegraphics[width=\linewidth]{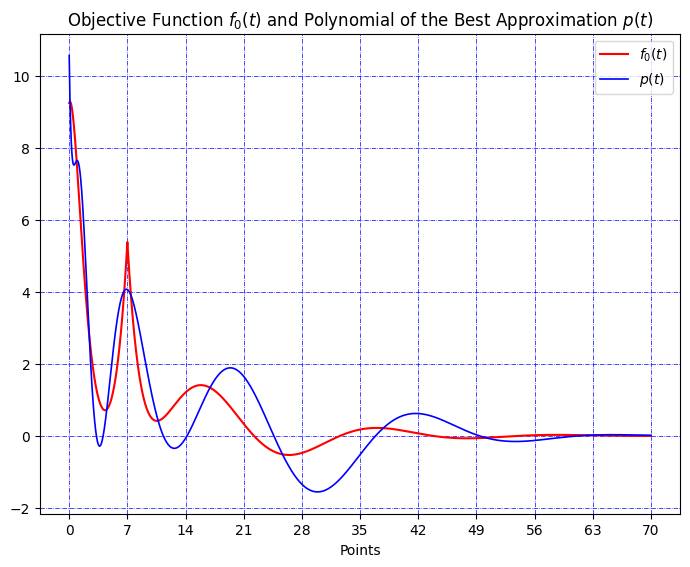}
			\caption{$ f(t) $																																}
	\end{subfigure}}\hfill
	\adjustbox{valign=m}{\begin{subfigure}{0.3\textwidth}
			\centering
			\includegraphics[width=\linewidth]{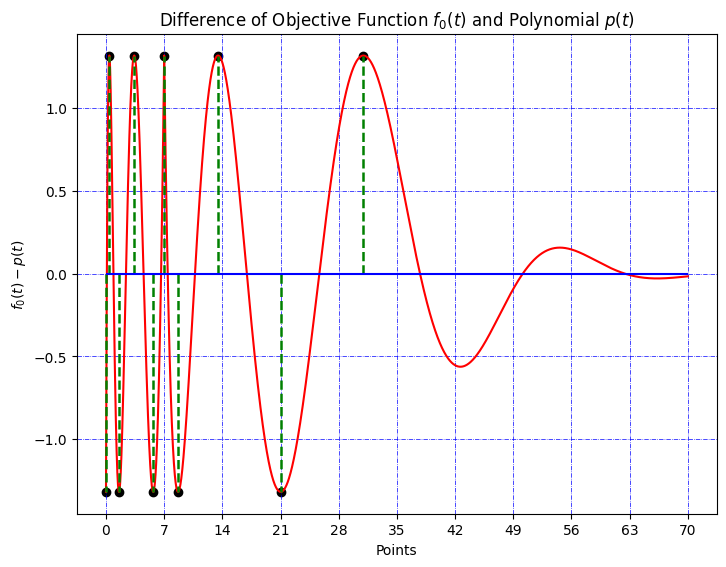}
			\caption{$ f(t) - p(t) $}
	\end{subfigure}}\hfill
	\adjustbox{valign=m}{\begin{subfigure}{0.3\textwidth}
			\centering
			\includegraphics[width=\linewidth]{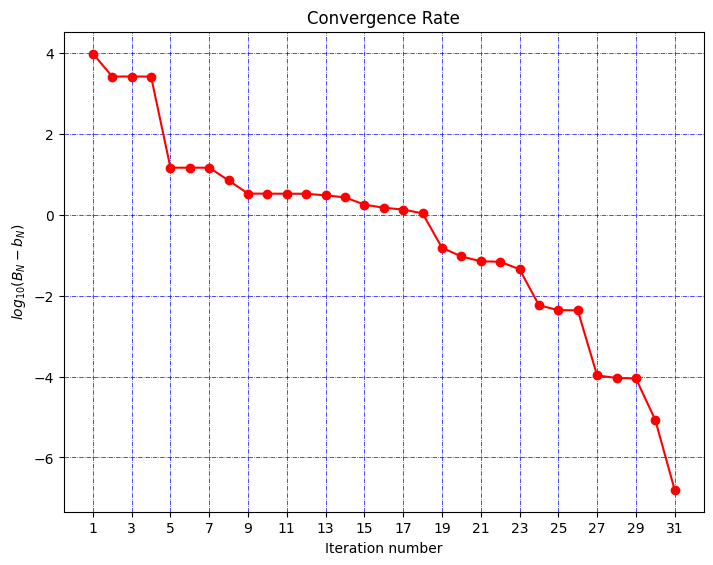}
			\caption{$ \log_{10} (B_N - b_N) $}
	\end{subfigure}}
	\label{fig:approximation_exponents}
\end{figure}

\textbf{Unconstrained problems.} We consider the space~$\cP$ of 
exponential polynomials on~$[0,+\infty)$: 
$$
p(t) \quad = \quad  
e^{-0.1t} \bigl( a_3 \,\cos \, (0.2 t) \, + \, a_4 \, \sin \,(0.2 t)  \, + 
a_5 \, \cos \, (0.3 t) \, + \, a_6 \, \sin \,(0.3 t)\, \bigr)  
$$
$$
+\ e^{-0.5t} \, \bigl( a_1 \, \cos \, (0.4 t) \, + \, a_2 \, \sin  \, (0.4 t) \bigr)  \ 
+ \  e^{-0.9t} \bigl( a_7 \, \cos \, t \, + \, a_8 \, \sin t\, \bigr) 
\ 
+ \  a_9\, e^{-0.3t}
$$ 
and a signal~$f$, which is a sum of these exponentials with coefficients: 
$$
(a_k)_{k=1}^9 \ = \ \Bigl(\, 1\, , \,  1\, , \, 
4\, , \,  -7\, , -3\, , \, -2\, , \,  1\, , \,  5\, , \,  6 \Bigr)
$$
and of the noise 
$\, 8\, e^{\frac{-|t-7|}{2}}$. Algorithm~2 finds the best approximating polynomial~$p \in \cP$ 
with the precision~$\varepsilon = 10^{-8}$
within~$k=31$ iterations,  
Fig.~\ref{fig:approximation_exponents}. 
The distance~$\|p-f\|_{C(\re_+)}\,  = \, 1.318352 $ (we leave six digits after the decimal point). The algorithm converges slower than in the previous example (subsection 9.2) because of the unbounded domain and a slow decrease of the 
functions~$\varphi_1, \ldots , \varphi_4$ (with the exponential part~$e^{-0.1t}$). 

\textbf{Problems with linear constraints.}
The same problem under the constraint $\, \int_{\re_+} p\, dt\,  \, = \, 1$
gives the following results. 
The distance to the best approximation~$\, \|p-f\|\,$  increases in more than 
$1.5$ times, 
from~$\, 1.318352 \, $ to $ \, 2.104564$.
The best approximation polynomial~$p$ now possess a 
degenerated $L$-alternance that    
 consists of $5$ points instead of~$n+1-r = 9$.  
Because of the degeneracy, Algorithm~1c stacks on this problem and we have to 
apply the regularization. Algorithm~2c finds the solution  within~$43$ iterations 
with the precision~$\varepsilon = 10^{-6}$.    

\bigskip

\begin{center}
	\textbf{9.3. Recovery of highly non-stationary signal}
\end{center}
\bigskip

The following  signal on the segment~$[0,1]$ was analysed 
in~\cite{Huang}:   
\begin{equation}
	\label{non-stationary signal}
	f(t) \quad = \quad \underbrace{\cos(4 \pi \lambda(t) t)}_{\varphi_1(t)} \quad + \quad 2\, \underbrace{\sin(4 \pi t)}_{\varphi_2(t)}\, ,
\end{equation}
where 
\begin{equation*}
	\lambda(t) = 
	\begin{cases}
		4 \ + \ 32 t \qquad  \ \ , \quad 0   \leqslant t \leqslant 0.5, \\
		4 \ + \ 32(1 - t)\ , \quad 0.5 < t \leqslant 1.
	\end{cases}
\end{equation*}

The polynomial~$p$ of the best approximation by the system~$\Phi = \{\varphi_1, \varphi_2\}$
is, of course, equal to~$f$, i.e.,~$\|p-f\|_{C[0, 1]} = 0$. 
Algorithm~2 finds this solution in two iterations.

\begin{figure}[h]
	\caption{}
	\centering
	\adjustbox{valign=m}{\begin{subfigure}{0.3\textwidth}
			\centering
			\includegraphics[width=\linewidth]{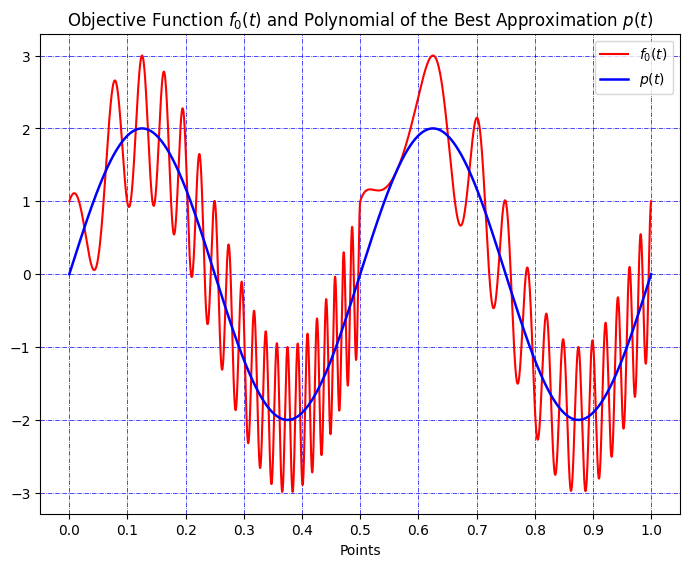}
			\caption{$ f(t) $ and $ p(t) $}
	\end{subfigure}}\hfill
	\adjustbox{valign=m}{\begin{subfigure}{0.3\textwidth}
			\centering
			\includegraphics[width=\linewidth]{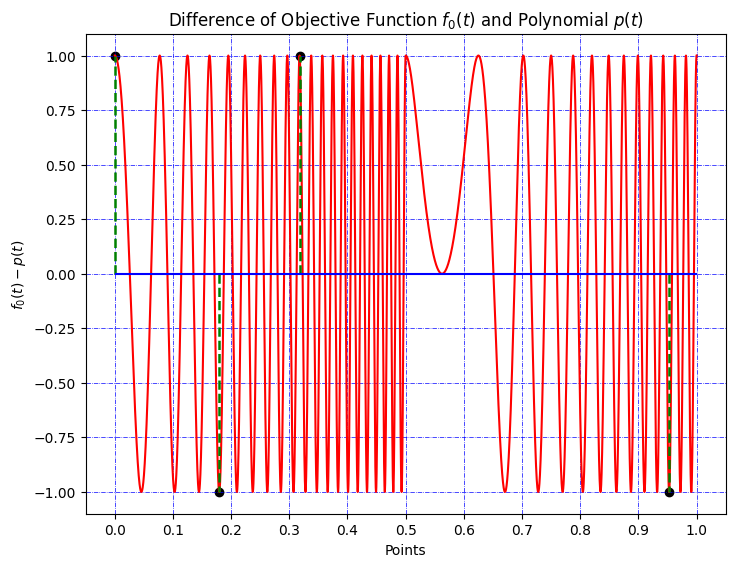}
			\caption{$ f(t) - p(t) $}
	\end{subfigure}}\hfill
	\adjustbox{valign=m}{\begin{subfigure}{0.3\textwidth}
			\centering
			\includegraphics[width=\linewidth]{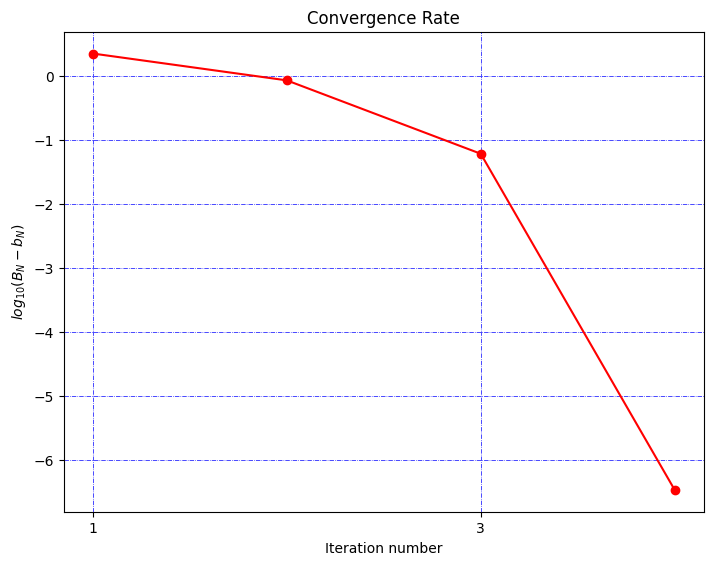}
			\caption{$ \log_{10} (B_N - b_N) $}
	\end{subfigure}}
	\label{fig:signal_from_paper_with_EMD}
\end{figure}

Approximating the same signal by the system~$\phi_1 = 1, \phi_2 = \cos(4 \pi t) , \phi_3 = \sin(4 \pi t) $ is shown in Fig.~\ref{fig:signal_from_paper_with_EMD}. The best approximation polynomial is $ p(t) \ = \  2 \sin(4 \pi t) $, Algorithm~2 finds it within three iterations. 
In fact, it  decomposes the signal $ f(t) $ into the sum of  $  2 \sin(4 \pi t) $ 
and $ f(t) - 2 \sin(4 \pi t) = \varphi_1(t) $ as the Empirical Mode Decomposition (EMD)  using a non-uniform mask~\cite[Sec. Numerical Experiments]{Lin_EMD}. 
\bigskip

\begin{center}
	\textbf{9.4. Markov-Bernstein type inequalities}
\end{center}
\bigskip 

We apply Algorithm~2c  to find the sharp constants in the Markov-Bernstein inequalities
(subsection~8.2) for various functional 
systems. 

\smallskip 

\textbf{Lacunary algebraic polynomials}  by the system of power functions
$\Phi \, = \,  \{t^{m_k}\}_{k=1}^n$. We apply Algorithm~2c to the problem~(\ref{eq.MB}) 
with~$ \varepsilon = 10^{-6}$. By Proposition~\ref{p.30} the constant $C_j$ 
is the reciprocal to the solution of this problem.

\renewcommand{\arraystretch}{2.1}

\begin{table}[h]
	\centering
	\caption{The Markov-Bernstein constants $ C_1, C_2 $ for different systems of lacunary algebraic polynomials.}
	\begin{tabular}{|c|c|c|}
		\hline 
		Basis of lacunary polynomials & $ C_1 $ & $ C_2 $ \\
		\hline
		$1, t, t^2, t^3, t^4, t^5, t^6 $  & $ 36 $ & $ 420  $ \\
		\hline
		$ 1, t, t^2, t^3, t^5, t^6 $  & $ 25.060144  $ & $ 201.979398$ \\
		\hline
		$ 1, t, t^{3}, t^5, t^6 $  & $ 25  $ & $ 200 $ \\
		\hline
		$ 1, t, t^5, t^6 $  & $  13.831259  $ & $ 69.1085 $ \\
		\hline
		$ 1, t, t^6 $  & $ 12  $ & $ 60 $ \\
		\hline
	\end{tabular}
	\label{table:Markov_Bernstein}
\end{table}

Note that if we remove monomials of even degree from the Chebyshev system $ \mathcal{C}_6 = \{1, t, t^2, t^3, t^4, t^5, t^6\} $, the constants $ C_1, C_2 $ become closer to the corresponding constants for the Chebyshev system $ \mathcal{C}_5 = \{1, t, t^2, t^3, t^4, t^5 \}$ (Table \ref{table:Markov_Bernstein}). 

\bigskip 

\textbf{Complex exponentials}. 
We find  the Markov-Bernstein constants for the space of 
quasipolynomials~$\cP \, = \, 
\bigl\{\, e^{-t} \bigl( a_1 \, \cos\, t\, + \, + a_2\,\sin\, t \, + \, 
a_3  \,  \quad a_1, a_2, a_3 \in \re \, \bigr\}$. Algorithm~2c applied to the problem 
~$\|p\|_{C(\re_+)} \to \min\, , \ p'(0) = 1$
returns the coefficients  
$a_1 = 1.006772 ., \, a_2 \, = \, 0.884983,  \, a_3 \, = \, 1.121789$. By~Proposition~\ref{p.40}, the corresponding Markov-Bernstein constant is~$ C_1 = \dfrac{1}{\|p(t)\|_{C(\re_+)}} = 8.694367 $. The alternance is nondegenerate and consists of~$m=3$ points, the precision $\varepsilon = 10^{-6}$ is achieved within~$8$ iterations of Algorithm~2c.
The algorithm runs $1.7$ sec. in a standard laptop.

\bigskip

\begin{center}
	\textbf{9.5. How often does the degeneracy occur?}
\end{center}
\bigskip

In our numerical experiments with polynomials by complex exponentials, shifts of Gaussian and Cauchy functions, 
lacunary algebraic and trigonometric polynomials,  
in~12 $\%$ of non-constrained problems degenerate simplices occur in Algorithm~1, 
after which the algorithm suffers and we had to apply the regularization and Algorithm~2. 
Moreover,~8 $\%$ of those cases have degenerate alternances, which consists of less than $n+1$ points. For problems with linear constraints, there number of degenerate cases is much larger:    
in 82~$\%$ of cases degenerate simplices  appear in Algorithm~2c,  and 73~$\%$ of cases produce  degenerate alternances. 
 
In this section we show the statistics over random system of functions, and it turns out to be quite similar.  
For this purpose, we generate $ m $ random independent uniformly distributed points 
(knots) on~$[-1,1]$, order them as $-1< t_1 < \cdots < t_m < 1$, then 
generate other $m$ random points $v_1, v_2, \ldots v_m $ (the values of the function) on the same segment $ [-1, 1] $,    interpolate the points $ (t_1, v_1), \ldots, (t_m, v_m) $ with a cubic  spline $ \varphi_1(t) $. This is a unique interpolating~$C^2[-1, 1]$ spline 
with the  ``not-a-knot'' boundary conditions, i.e., 
the two extreme polynomials coincide on each side of the segment. 

In the same way 
we generate functions~$\varphi_2 , \ldots , \varphi_n$.  
The boundary conditions are chosen so that the first and second segment at a curve end are the same polynomial. 
Having generated the system~$\Phi$ on the segment~$[-1,1]$, we test Algorithm~2 with 
finding the best approximations for  two functions: 1) $f(t) = |t| $ and 2) a random function~$f$
generated as above (in the same way as all~$\varphi_i$). Then we test Algorithm~2c with the  problem $\|p\|_{C[-1,1]} \to \min$ under the constraint~$\ell(p) = \sum_{i=1}^n\, p^{i}\, = \, 1$. 
For each pair $(m, n) $ we run 100 random systems~$\Phi$ 
and then calculate the proportion of degenerate cases, average number of iterations and average program runtime for degenerate and non-degenerate cases. We call a case degenerate if at least in one iteration we have~$ q_{i_1, i_2}(t_0) \, < \, 0.05 $ for some~$i_1\ne i_2$, where 
$t_0$ is the point of maximum of~$|p_k(t) - f(t)|$ on the segment~$[-1,1]$. 
The polynomial~$q_{i_1, i_2}$ is normalized as~$\|\bq_{i_1, i_2}\| = 1$. The results  are given in Table \ref{table:Cubic_Splines}.

\renewcommand{\arraystretch}{1}

\begin{table}[h]
	\footnotesize
	\centering
	\caption{The performance of Algorithms 2 and 2c for random polynomial systems}
	\begin{tabular}{|c|c|c|c|c|c|c|c|}
		\hline 
		\multirow{2}{*}{Problem} & \multirow{2}{*}{Parameters} & \multicolumn{3}{c|}{Non-degenerate cases}  &  \multicolumn{3}{c|}{Degenerate cases} \\
		\cline{3-8}
		& & share  & iter. & time (sec) & share & iter. &  time (sec) \\
		\hline
		\multirow{3}{*}{$f=0, \, \ell(p) = 1$} &
		$ m=10, n=3 $  & $ 0.83 $ & $ 5.66 $ & $ 1.75 $ &  $ 0.17 $ & $ 19.41 $ & $ 6.85 $ \\
		\cline{2-8}
		& $ m=10, n=5 $  & $ 0.59 $ & $ 12.02 $ & $ 4.14 $ &  $ 0.41 $ & $ 18.95$ & $ 7.26 $ \\
		\cline{2-8}
		& $ m=5, n=7 $  & $ 0.27 $ & $ 11.04 $ & $ 3.75 $ &  $ 0.73 $ & $ 24.79 $ & $ 10.24 $ \\
		\hline
		
		\multirow{3}{*}{$ f = |t| $} &
		$ m=10, n=3 $  & $ 0.33 $ & $  10.81 $ & $ 8.27 $ &  $ 0.67 $ & $ 17.27 $ & $ 14.51 $ \\
		\cline{2-8}
		& $ m=10, n=5 $  & $ 0.17 $ & $ 15.5 $ & $ 10.38 $ &  $ 0.83 $ & $ 22.12 $ & $ 21.79 $ \\
		\cline{2-8}
		& $ m=5, n=7 $  & $ 0.37 $ & $  11.47 $ & $ 10.56 $ &  $ 0.63 $ & $ 17.29 $ & $ 20.18 $ \\
		\hline
		
		\multirow{3}{*}{$f$ is random} &
		$ m=10, n=3 $  & $ 0.63 $ & $ 8.16 $ & $ 5.51 $ &  $ 0.37 $ & $ 19.24 $ & $ 16.21 $ \\
		\cline{2-8}
		& $ m=10, n=5 $  & $ 0.36 $ & $ 13.74 $ & $ 8.13 $ &  $ 0.64 $ & $ 26 $ & $ 17.85 $ \\
		\cline{2-8}
		& $ m=5, n=7 $  & $ 0.39 $ & $ 12.47 $ & $ 6.76 $ &  $ 0.61 $ & $ 24.31 $ & $ 15.39 $ \\
		\hline
	\end{tabular}
	\label{table:Cubic_Splines}
\end{table}

Several empirical conclusions can be made from Table \ref{table:Cubic_Splines}. 
The share of degenerate cases
increases in the degree~$n$ of the system. Also approximation of a non-smooth function 
(in our case~$f=|t|$) admits more  degenerate cases. 
The number of iterations in degenerate cases is, on average, from 1.5 to 2 times greater than in non-degenerate ones.
\newpage

\begin{center}
	\large{\textbf{10. Discussion. Recursive approximations vs convex programming}}
\end{center}
\bigskip

Both constrained and unconstrained problems of uniform polynomial approximation 
can be solved by tools of convex optimisation. However, in practice they are 
less efficient that Remez-like recursive methods. In particular, for 
approximating by Chebyshev systems, in particular, for algebraic polynomial approximation,  everybody uses the Remez algorithm.  The main reason is that 
it gives an a posteriori estimate and hence, usually reaches a desired accuracy 
in much fewer iterations  than it is prescribed theoretically. 
Another reason is that convex optimization algorithms are 
usually suffer because of the non-smoothness of this problem. 
Actually, it belongs to the class of minimax problems: 
\begin{equation}
	\label{min_max_formulation}
	\min_{p_1, \ldots, p_n \in  \R} \max_{t \in \mathcal{K}}| f(t) - \sum_{i=1}^n p_i \varphi_i(t) |.
\end{equation}
A popular way to solve it by the  gradient descent ascent method (GDA) with the 
iteration: 
\begin{equation*}
	x_{k+1} = x_k - \alpha_k \nabla_x F(x_k, y_k), \quad y_{k+1} = y_k + \beta_k \nabla_y F(x_k, y_k); \quad \alpha_k, \beta_k \in \R_{\geqslant 0}, 
\end{equation*}
Usually it is
quite efficient for minimax problems~\cite{Fallah, Liang, Zhang}.  
However, in our case  the objective function $ F(p, t) = |f(t) - \sum_{i=1}^n p_i \varphi_i(t) | $ is non-smooth and not necessarily concave with respect to variable $ t $. Most existing literature deals with GDA algorithms for solving minimax problems of the form $ \min_{x \in \mathcal{X}} \max_{y \in \mathcal{Y}} F(x, y) $, where the objective function $ F(x, y) $ is at least smooth and  $ F(x, \cdot) $ is concave for each $ x \in \mathcal{X} $ \cite{Lin, Nouehed}.

\bigskip

\noindent {\large \textbf{Acknowledgements}}
\smallskip 

The research  is performed with the support of the Theoretical Physics and Mathematics Advancement Foundation ``BASIS'', grant 22-7-1-20-1. 
\bigskip 

\noindent  {\large \textbf{Declarations}}
\smallskip

\textbf{Conflict of interest}. The authors have no interests to delcare.

\end{document}